\documentclass{scrartcl}

\KOMAoptions{paper=a4}
\KOMAoptions{fontsize=12pt}
\KOMAoptions{DIV=calc} 

\usepackage[utf8]{inputenc} 
\usepackage[T1]{fontenc} 
\usepackage[ngerman,english]{babel} 
\usepackage{csquotes} 
\usepackage{hyphenat} 

\usepackage{kpfonts}

\usepackage{cabin}

\usepackage{microtype} 

\addtokomafont{title}{\rmfamily}

\usepackage{enumitem}
\setlist[enumerate]{label*=(\alph*),ref=(\alph*),itemsep=0pt,topsep=5pt}
\setlist[itemize]{itemsep=0pt,topsep=5pt}

\KOMAoptions{DIV=calc}

\usepackage{booktabs} 

\usepackage{graphicx} 
\graphicspath{{pic/}} 

\usepackage[format=plain,labelfont={sf,footnotesize},textfont=small,margin=12pt]{caption} 
\DeclareCaptionLabelSeparator{MySpace}{\enskip } 
\usepackage[style=alphabetic,
						sorting=nyt,
						maxnames=4,
						backend=biber,
						doi=true,
						url=true,
						isbn=false,
						clearlang=true]
						{biblatex}

\bibliography{C:/Users/j.rau/Dropbox/Bibliography/biblio}

\DeclareFieldFormat*{title}{\textit{#1}} 
\renewbibmacro{in:}{} 
\DeclareFieldFormat{journaltitle}{\iffieldequalstr{journaltitle}{ArXiv e-prints}{Preprint}{#1}} 
\DeclareFieldFormat{booktitle}{#1} 
\AtEveryBibitem{\clearfield{eprintclass}} 
\DeclareRedundantLanguages{english,English}{english,german,ngerman,french}
\DeclareSourcemap{ 
  \maps[datatype=bibtex]{
    \map{ 
      \step[fieldset=addendum,null]
      \step[fieldset=eprintclass,null]
    }  
    \map{ 
      \step[fieldsource=doi,final]
      \step[fieldset=url,null]
    }  
  }
}

\usepackage{hyperref} 
\usepackage{xcolor} 
\hypersetup{breaklinks=true} 
\definecolor{darkblue}{RGB}{0,0,170}
\definecolor{darkred}{RGB}{200,0,0}
\hypersetup{citecolor=darkblue, urlcolor=darkblue, linkcolor=darkblue, colorlinks=true} 
\usepackage[all]{hypcap} 

\usepackage{aliascnt} 
\usepackage[amsmath,thmmarks]{ntheorem} 

\newtheorem {theorem}{Theorem}[section]

\newaliascnt{proposition}{theorem}
\newtheorem {proposition}[proposition]{Proposition}
\aliascntresetthe{proposition}

\newaliascnt{lemma}{theorem}
\newtheorem {lemma}[lemma]{Lemma}
\aliascntresetthe{lemma}

\newaliascnt{corollary}{theorem}
\newtheorem {corollary}[corollary]{Corollary}
\aliascntresetthe{corollary}

\newaliascnt{conjecture}{theorem}

\aliascntresetthe{conjecture}

\theorembodyfont{\upshape}

\newaliascnt{definition}{theorem}
\newtheorem {definition}[definition]{Definition}
\aliascntresetthe{definition}

\newaliascnt{example}{theorem}
\newtheorem {example}[example]{Example}
\aliascntresetthe{example}

\newaliascnt{exercise}{theorem}

\aliascntresetthe{exercise}

\newaliascnt{goal}{theorem}

\aliascntresetthe{goal}

\newaliascnt{construction}{theorem}

\aliascntresetthe{construction}

\theoremheaderfont{\itshape}
\theorembodyfont{\upshape}

\newaliascnt{remark}{theorem}
\newtheorem {remark}[remark]{Remark}
\aliascntresetthe{remark}

\newaliascnt{convention}{theorem}

\aliascntresetthe{convention}

\newaliascnt{notation}{theorem}

\aliascntresetthe{notation}

\theoremstyle {nonumberplain}
\theoremheaderfont{\itshape}
\theorembodyfont{\upshape}
\theoremsymbol{\ensuremath{_\blacksquare}}

\newtheorem {proof}{Proof}

\addto\extrasenglish{%
}

\usepackage {amsmath} 
\usepackage{tikz}
\usetikzlibrary{cd}
\usetikzlibrary{babel}

\newcommand{\C}{{\mathbf C}}

\renewcommand{\P}{{\mathbf P}}
\newcommand{\Q}{{\mathbf Q}}
\newcommand{\R}{{\mathbf R}}
\newcommand{\T}{{\mathbf T}}

\newcommand{\Z}{{\mathbf Z}}

\newcommand{\TP}{\mathbf{TP}}

\newcommand{\DD}{{\mathcal D}}
\newcommand{\EE}{{\mathcal E}}

\newcommand{\HH}{{\mathcal H}}

\newcommand{\OO}{{\mathcal O}}

\renewcommand{\SS}{{\mathcal S}}

\newcommand{\Bl}{{\mathcal Bl}}

\newcommand{\wt}[1]{{\widetilde{#1}}}

\newcommand{\wtQ}{{\widetilde{Q}}}

\newcommand{\BlX}{{\widetilde{X}}}
\newcommand{\Blpi}{{\widetilde{\pi}}}
\newcommand{\BlI}{{\widetilde{I}}}
\newcommand{\BlDD}{{\widetilde{\DD}}}
\newcommand{\Rneg}{\R_{\leq 0}}
\newcommand{\Rpos}{\R_{\geq 0}}
\newcommand{\cd}{k}

\DeclareMathOperator{\Star}{Star}

\DeclareMathOperator{\codim}{codim}

\DeclareMathOperator{\Rec}{Rec}
\DeclareMathOperator{\Hom}{Hom}
\DeclareMathOperator{\Trop}{Trop}

\makeatletter
\def\widebreve{\mathpalette\wide@breve}
\def\wide@breve#1#2{\sbox\z@{$#1#2$}%
     \mathop{\vbox{\m@th\ialign{##\crcr
\kern0.08em\brevefill#1{0.8\wd\z@}\crcr\noalign{\nointerlineskip}%
                    $\hss#1#2\hss$\crcr}}}\limits}
\def\brevefill#1#2{$\m@th\sbox\tw@{$#1($}%
  \hss\resizebox{#2}{\wd\tw@}{\rotatebox[origin=c]{90}{\upshape(}}\hss$}
\makeatletter

\addtokomafont{author}{\large}
\addtokomafont{section}{\large}
\addtokomafont{subsection}{\normalsize}

\begin {document}

\title {Real semi-stable degenerations, real-oriented blow-ups and straightening corners}
\author {Johannes Rau}
\date{} 

\maketitle

\begin{abstract}
	\noindent
  We study totally real semi-stable degenerations (and more generally, smooth semi-stable degenerations). 
	Our goal is to describe	the homeomorphism type of the real locus $\R X_t$ of the general fibre in terms of the special fibre. 
	We give a general homeomorphism statement via the real-oriented blow-up of the family. 
	Using this, we give more explicit descriptions of $\R X_t$ as a stratified space glued from 
	(covers of) strata of the special fibre. We also give relative versions of the statements
	and consider the example of toric degenerations in order to link the technique to tropicalisation. 
\end{abstract}


\section{Introduction}

Let $\pi \colon X \to \DD$ be a complex semi-stable degeneration 
over the unit disc $\DD \subset \C$.  
Assume that $\pi$ is real, that is, $X$ is equipped with an anti-holomorphic involution 
compatible with complex conjugation on $\DD$ via $\pi$. 
We say $\pi$ is totally real if locally each branch of the normal crossing divisor
over $0$ is real, that is, fixed by the real structure. 
The detailed definitions can be found in the main text. 

\begin{theorem} \label{thm:MainComplex}
  Let $\pi \colon X \to \DD$ be a totally real semi-stable degeneration and let  
	$\Blpi \colon \BlX \to \BlDD = S^1 \times [0,1)$ be the real-oriented blow-up of $\pi$. 
	We set $\BlX^+ := \Blpi^{-1}(\{+1\} \times [0,1))$ and $\BlX_0^+ := \Blpi^{-1}(\{+1\} \times \{0\})$. 
	Then there exists a homeomorphism
	\begin{equation} 
		H \colon \R \BlX_0^+ \times [0,1) \to \R \BlX^+
	\end{equation}
	such that $\Blpi(H(x,t)) = t$ for all $(x,t) \in \BlX_0^+ \times [0,1)$. 	
	In particular, the real positive special fibre $\R\BlX_0^+$ is homeomorphic 
	to the real locus $\R X_t$ of the generic real fibre $X_t$, $t \in (0,1)$.
\end{theorem}

See \cite[Theorem 5.1]{NO-RelativeRoundingToric} and 
\cite[Proposition 6.4]{Arg-RealLociLog} for similar, more general, 
statements in the setting of log smooth maps and the Kato-Nakayama 
space (see below for more comments). 
We also prove the following relative version of this theorem.

\begin{theorem} \label{thm:RelativeReal}
  \sloppy Let $\pi \colon X \to \DD$ be a totally real semi-stable degeneration and 
	$X_1, \dots, X_m \subset X$ a transversal collection of real submanifolds. 
	Then there exists a homeomorphism
	$H \colon \R \BlX_0^+ \times [0,1) \to \R \BlX^+$ as in \autoref{thm:MainComplex}
	such that 
	\begin{equation} 
		H\left(\R\BlX^+_{i,0} \times [0,1)\right) = \R\BlX^+_{i}
	\end{equation}
	for all $i = 1, \dots, m$. 
	In particular, 
	the topological tuples $\R\BlX^+_{1,0}, \dots, \R\BlX^+_{m,0} \subset \R\BlX^+_0$ and 
	$\R X_{1,t}, \dots, \R X_{m,t} \subset \R X_t$ are homeomorphic for all $t \in (0,1)$. 
\end{theorem}

In view of the above theorems, it is desirable to obtain concrete descriptions
of the real positive special fibre $\R \BlX_0^+$. Note that in this text, 
strata of a stratification are always understood to be pairwise disjoint, locally
closed subsets. 

\begin{theorem} \label{thm:StrataComplex}
	The real positive special fibre $\R \BlX_0^+$ of a totally real semi-stable degeneration 
	$\pi \colon X \to \DD$ is a Whitney-stratified space whose strata
	are (the connected components of) topological covers of degree $2^{k-1}$ of
	each stratum of $\R X_0$ of codimension $k = k(S)$
	with respect to $\R X$. 
	In particular, 
	\begin{equation} 
		\chi^c(\R \BlX^+_0) = \sum_{k=1}^n \sum_{\substack{S \in \SS \\ k(S) = k}} 2^{k-1} \chi^c(S).
	\end{equation}  
	Here, $\SS$ is any refinement of the canonical 
	normal crossing Whitney stratification $\SS(\R X, \R X_0)$ on 
	$\R X$, but $k(S)$ always denotes the codimension
	with respect to $\SS(\R X, \R X_0)$
	(see \autoref{secStratifying}, page \pageref{secStratifying}, for detailed definitions). 
	In particular, if $X_1, \dots, X_m \subset X$ is a transversal collection of real divisors,
	we can choose $\SS$ to be the stratification induced by the normal crossing divisor 
	$\R X_0 \cup \R X_1 \cup \dots \cup \R X_m$. 
\end{theorem}

A more explicit recipe for how the strata are glued together in $\R \BlX^+_0$ can be found
in \autoref{describing} 
(especially for the case when all strata $S \in \SS$ are simply connected). 

\paragraph{Motivation and context}
The idea of using simple (in this case, semi-stable) degenerations to study
the generic fibre is of course very old and well-known. 
In particular, it is present from the very beginning
in the topological study of real algebraic varieties. 
Axel Harnack's and David Hilbert's original \enquote{small perturbation} method
\cite{Har-UeberDieVieltheiligkeit,Hil-UeberDieReellen}
is based on generic pencils of curves and hence semi-stable degenerations,
see \autoref{ex:smallperturbation}. 
A much more powerful incarnation of that idea was then introduced
for hypersurfaces in toric varieties
in the form of Oleg Viros's patchworking method, 
see for example \cite{Vir-CurvesDegree7, Vir-PatchworkingRealAlgebraic}.
Again, semi-stable degenerations appeared at least implicitly in the proofs,
e.g.\ \cite[Theorems 2.5.C and 2.5.D]{Vir-PatchworkingRealAlgebraic} and \autoref{ex:patchworking}.
We recommend \cite{Vir-IntroductionTopologyReal} for a beautiful account of
these ideas and their history. 
Other types of semi-stable deformations, in particular, 
deformations to the normal cone, see \autoref{ex:normalcone}, 
are also commonly used in real algebraic geometry,
often cleverly combined with Viro's patchworking method, e.g.\ 
\cite{Shu-LowerDeformationsIsolated,ST-PatchworkingSingularAlgebraic,BDIM-RealAlgebraicCurves}.
Inspired by these results, in this paper we search for a 
\enquote{patchwork} description of the
generic fibre for general semi-stable degenerations. 

The connection to Viro's patchworking 
will be made even more explicit in our upcoming joint work 
\cite{RRS-RealPhaseStructuresb} 
with Arthur Renaudineau and Kris Shaw.
Based on our previous work \cite{RRS-RealPhaseStructures}, 
we plan to apply the techniques of this paper to tropicalisations
of families of real algebraic varieties $X$. 
Inspired by the 
tropical description of Viro's patchworking
(cf.\ \cite{Vir-DequantizationRealAlgebraic, Mik-AmoebasAlgebraicVarietiesb, IMS-TropicalAlgebraicGeometry}),
the goal is to describe 
the homeomorphism type of the real locus of the generic fibre 
of $X$ in terms of its associated tropical variety $\Trop(X)$ 
equipped with a real phase structure $\EE(X)$, 
see \cite[Definition 2.2]{RRS-RealPhaseStructures}.
Under suitable smoothness assumptions, the situation can be reduced to 
semi-stable toric degenerations and hence to the statements in this paper,
see \autoref{cor:tropicalisation}. 
The idea to use the approach via real semi-stable degenerations was brought to our attention
by Erwan Brugallé in his recent work \cite{Bru-EulerCharacteristicSignature}. 
In this work, Brugallé uses the technique to prove that 
real algebraic varieties close to smooth tropical limits satisfy 
$\chi = \sigma$ (the Euler characteristic of the real part is equal to the signature of the complex part),
generalising previous works on patchworks in  
\cite{Ite-TopologyRealAlgebraic,Ber-EulerCharacteristicPrimitive,BB-EulerCharacteristicReal}.
With \autoref{thm:RelativeReal} and \autoref{thm:StrataComplex} 
we provide the slightly refined ingredient for the argument 
on the Euler characteristic side. 

Finally, one should mention that the techniques used here
are in fact very similar to a significantly more general
construction, studied in the language of logarithmic geometry, of 
the so-called \emph{Kato-Nakayama space} or \emph{Betti realization}
$X_{\text{log}}$
associated to a log analytic space $X$, see 
\cite[(1.2)]{KN-LogBettiCohomology},
also \cite{NO-RelativeRoundingToric}.
Indeed, if $X$ is the log space associated
to a semi-stable degeneration, its
Kato-Nakayama space $X_{\text{log}}$
coincides with the real-oriented blow-up,
see \cite[(1.2.3)]{KN-LogBettiCohomology}.
The study of the real locus in presence of a real structure
on $X$ is done
in \cite{Arg-RealLociLog}.
In the overlap of the settings, \autoref{thm:MainComplex}
is a special incarnation of 
\cite[Theorem 5.1]{NO-RelativeRoundingToric} and 
\cite[Proposition 6.4]{Arg-RealLociLog}.

\subsection*{Acknowledgement}

My first thanks go to Arthur Renaudineau and Kris Shaw for
the previous and ongoing collaborations which led me to this side project. 
Many joint conversations and discussions have inspired this text. 
Many thanks go to Erwan Brugallé for being the second source of inspiration, 
for many useful discussions on the subject, for corrections on earlier versions 
of the text and for sharing many ideas that entered here, in particular, 
\autoref{ex:normalcone} and \autoref{ex:deformationconic}.
It is my special pleasure to thank Stefan Behrens for useful discussions and comments 
and for pointing me to the references and techniques that are crucial in the technical core
of the paper in \autoref{straightening}. I would also like to thank Lionel Lang 
for discussions on real-oriented blow-ups related to a different project,
and Jean-Baptiste Campesato for mentioning to me early appearences
of real-oriented blow-ups and his related works. 
Finally, I would like to thank Helge Ruddat and Hülya Argüz
for pointing out to me the connections to logarithmic geometry
and Kato-Nakayama spaces.

The author is supported by the FAPA project \enquote{Matroids in tropical geometry} 
from the Facultad de Ciencias, Universidad de los Andes, Colombia.

\section{Smooth semi-stable degenerations} \label{SmoothSemiStable}

Our main statements are formulated in the context of complex semi-stable degenerations
since the applications we have in mind are situated in this context. 
However, to a large extent our story takes places solely in the smooth manifolds 
of real points inside these degenerations. 
We hence start by setting up the corresponding framework of \emph{smooth} 
semi-stable degenerations. 

Let $X$ be a smooth manifold. 
A subset $D \subset X$ is called a \emph{normal crossing divisor}
if for any point $p \in D$ there exists a chart $U \to V \subset \R^n$ centred in $p$
that identifies $D\cap U$ with $\{x_1 \cdots x_k = 0\} \cap V$ for some $1 \leq k \leq n$. 
We refer to such a chart as \emph{normal crossing chart of type $(k,n)$ for $D$ at $p$}.
Since any other chart at $p$ has the same type $(k,n)$, we may also speak of the type of $p$. 

A smooth function $\rho \colon X \to \R$ is called a \emph{normal crossing function}
if for any point $p$ with $\rho(p) = 0$ there exists a chart $U \to V \subset \R^n$
centred in $p$ such that $\rho(x) = x_1 \cdots x_k$ in local coordinates. 
Clearly, in such case $D = \rho^{-1}(0)$ is a normal crossing divisor. 
We refer to such a chart as \emph{normal crossing chart of type $(k,n)$ for $\rho$ at $p$}.

\begin{definition} \label{def:smoothSemiStable}
	A \emph{smooth semi-stable degeneration} is a proper normal crossing function  
	$\pi \colon X \to I$ 
	from a smooth manifold $X$ to the interval $I = (-1,1) \subset \R$ which is 
	submersive over $I \setminus \{0\}$. 
\end{definition}

We denote the fibre over $t \in I$ by $X_t = \pi^{-1}(t)$. 
A generic fibre $X_t$, $t \neq 0$, is a smooth submanifold of $X$. 
The special fibre $X_0$ is a normal crossing divisor in $X$. 
In fact, for most of what follows, it is sufficient
to assume that $X_0$ is normal crossing but drop the condition that $\pi$ 
is normal crossing, see \autoref{rem:weaklySemiStable}.
Our goal is to understand the homeomorphism type of a generic fibre 
via data contained in the special fibre $X_0$ (and the degeneration, of course). 
A systematic way of doing so is using the language of real-oriented blow-ups.

\begin{example} \label{ex:cubic}
  We use the following simple example to illustrate the upcoming definitions and statements.
	For more interesting examples, we refer to later sections. 
	We set $X = \R^2$ and consider the map $\pi \colon \R^2 \to \R$ given by
	\[
	  \pi(x,y) = y^2 - x^3 - 3x^2.
	\]
	Its two critical points are $(0,0)$ and $(-2,0)$.
	In the two-dimensional case, the existence of a $(2,2)$ chart
	at a critical point is clearly equivalent to the Hessian having 
	negative determinant. Hence $(0,0)$ admits such a chart while $(-2,0)$
	does not. Since	$\pi(-2,0) = -4$, it follows that 
	the restriction to $\pi^{-1}(I)$
	is a semi-stable degeneration, 
	see \autoref{semistablecubic}. 
	Conversely, restriction of the shifted function $\pi + 4$ 
	does not yield a semi-stable degeneration.
	\begin{figure}[ht]
		\begin{minipage}[r]{0.4\textwidth}
		  \centering
			\includegraphics[width=0.8\textwidth]{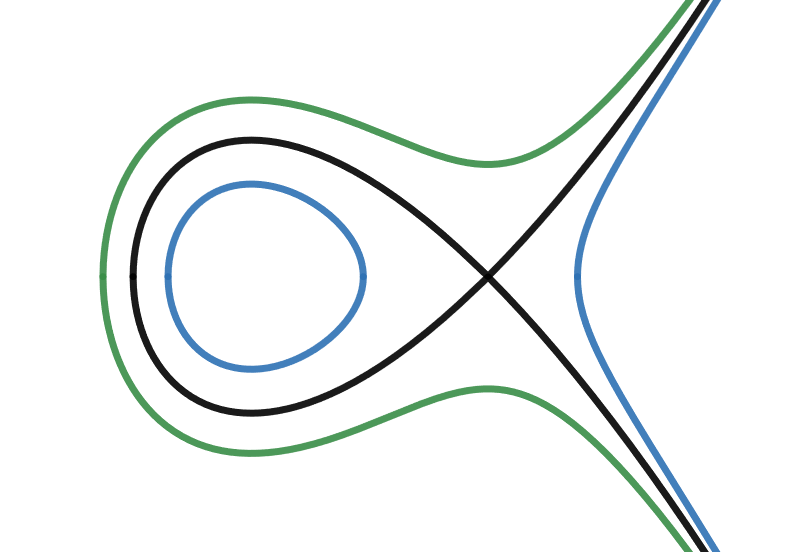}
		\end{minipage}\hfill
		\begin{minipage}[l]{0.5\textwidth}
			\caption{The special (black), positive (green) and negative (blue) 
			fibres of the smooth semi-stable degeneration $\pi(x,y) = y^2 - x^3 - 3x^2$}
			\label{semistablecubic}
		\end{minipage}		
	\end{figure}%
\end{example}

\section{Real-oriented blow-ups}

We start by reviewing the notion of real-oriented blow-up along a normal crossing divisor.
Our basic reference is \cite{HPV-CompactificationHenonMappings},
early references are \cite{AC-LaFonctionZeta,Per-DegenerationsAlgebraicSurfaces}. 

Let $X$ be a smooth manifold and $D \subset X$ a normal crossing divisor. 
We denote by $\Bl^\R (X, D)$
the \emph{real-oriented blow-up of $X$ along $D$} as defined in \cite[Definition 5.1, Proposition 5.2]{HPV-CompactificationHenonMappings}
(with the only difference that we are working in the smooth instead of real-analytic category here). 
For now, $\Bl^\R (X, D)$ is a topological space (later, we will equip it with the structure 
of a manifold with corners) together with a blow-down map $\alpha \colon \Bl^\R (X, D) \to X$
which is a homeomorphism over $X \setminus D$. 
It is unique up to canonical homeomorphism compatible with blow-down maps. 
In essence, $\Bl^\R (X, D)$ is a space that behaves like $X$ minus a \enquote{tubular neighbourhood} 
of $D$. 
Since the construction of the real-oriented blow-up is local, for our purposes 
it is enough to understand how it works on normal crossing charts. 
Note also that the following local description 
shows that $\alpha$ is proper. 
We set $\Rpos = [0, \infty)$ and $\Rneg = (-\infty, 0]$. 

\begin{lemma} \label{lem:localblowup}
  Set $X = \R^n$ and $D = \{x_1 \cdots x_k = 0\}$ for some $1 \leq k \leq n$. Then
	the real-oriented blow-up of $X$ along $D$ is given by
	\begin{equation} \label{eq:localblowup} 
	  (\Rneg \sqcup \Rpos)^k \times \R^{n-k} \to \R^k \times \R^{n-k}	
	\end{equation}
	where the map on the first $k$ factors is given by the obvious map $\Rneg \sqcup \Rpos \to \R$
	and on the last $n-k$ factors is identity. 
\end{lemma}

\begin{proof}
  This follow from \cite[Theorem 5.3]{HPV-CompactificationHenonMappings} and the base case
	\[
	  \Bl^\R (\R, \{0\}) = \Rneg \sqcup \Rpos \to \R
	\]
	discussed for example at \cite[Page 235]{HPV-CompactificationHenonMappings}.
\end{proof}

It will be convenient to label the $2^k$ orthants of $(\Rneg \sqcup \Rpos)^k$
in the following way: Given $\epsilon = (\epsilon_1, \dots, \epsilon_k) \in \{\pm\}^k$, we denote by $\R^k_\epsilon$ the (closed) orthant 
in $\R^k$ obtained as the closure of the set of points $x = (x_1, \dots, x_k) \in (\R^*)^k$
for which $x_i$ has sign $\epsilon_i$ for all $i = 1, \dots, k$.
Clearly,
\[
  (\Rneg \sqcup \Rpos)^k = \bigsqcup_{\epsilon \in \{\pm\}^k} \R^k_\epsilon.
\]
Whenever we want to explicitly specify the orthant that contains a point in $(\Rneg \sqcup \Rpos)^k$, 
we write this point as a tuple $(\epsilon, x)$, $x \in \R^k$ and $\epsilon \in \{\pm\}^k$.

The important fact for us is that real-oriented blow-ups are compatible with smooth
semi-stable degenerations, that is, we can blow-up the map $\pi$ in the following sense.

\begin{proposition} 
  Let $\pi \colon X \to I$ be a smooth semi-stable degeneration.
	Let $\alpha \colon \BlX = \Bl^\R (X, X_0) \to X$ and $\beta \colon \BlI = \Bl^\R(I,\{0\}) = (-1,0] \sqcup [0,1) \to I$ 
	denote the real-oriented blow-ups along the
	special fibre and along $0$, respectively. 
	Then there exists a unique continuous map $\Blpi \colon \BlX \to \BlI$
	such that the diagram 
	\[
	\begin{tikzcd}
		\BlX \arrow[r, "\alpha"] \arrow[d, "\Blpi"']
		& X \arrow[d, "\pi"] \\
		\BlI \arrow[r, "\beta"]
		& I
	\end{tikzcd}
	\]
	commutes. Moreover, $\Blpi$ is proper. 
\end{proposition}

\begin{proof}
  It suffices to prove the statement on $\pi$-normal crossing charts of type $(k,n)$.
	Hence it is enough 
	check the case $\pi : X = \R^n \to \R$, $\pi(x) = x_1 \cdots x_k$. 
	By \autoref{lem:localblowup}, the real-oriented blow-up is given by 
	$\alpha \colon (\Rneg \sqcup \Rpos)^k \times \R^{n-k} \to \R^k \times \R^{n-k}$. 
	Hence the only continuous map that agrees with $\beta^{-1} \circ \pi \circ \alpha$ 
	on $\alpha^{-1}( (\R^*)^k \times \R^{n-k})$ is given by
	\begin{equation} \begin{split} \label{eq:Blpi} 
		\Blpi \colon (\Rneg \sqcup \Rpos)^k \times \R^{n-k} &\to \Rneg \sqcup \Rpos, \\
		             ((\epsilon,x), y) & \mapsto (\epsilon_1 \cdots \epsilon_k, x_1 \cdots x_k)		
	\end{split} \end{equation}
	using the tuple notation from above. The properness follows easily from the fact that 
	$\alpha$ and $\pi$ are proper. 
\end{proof}

\begin{example} \label{ex:blowup}
  The real-oriented blow-up of the semi-stable degeneration $\pi(x,y) = y^2 - x^3 - 3x^2$ from \autoref{ex:cubic}
	is depicted in \autoref{blowup}. 
	
	\begin{figure}[h]
		\begin{minipage}[r]{0.4\textwidth}
		  \centering
			\includegraphics[width=0.8\textwidth]{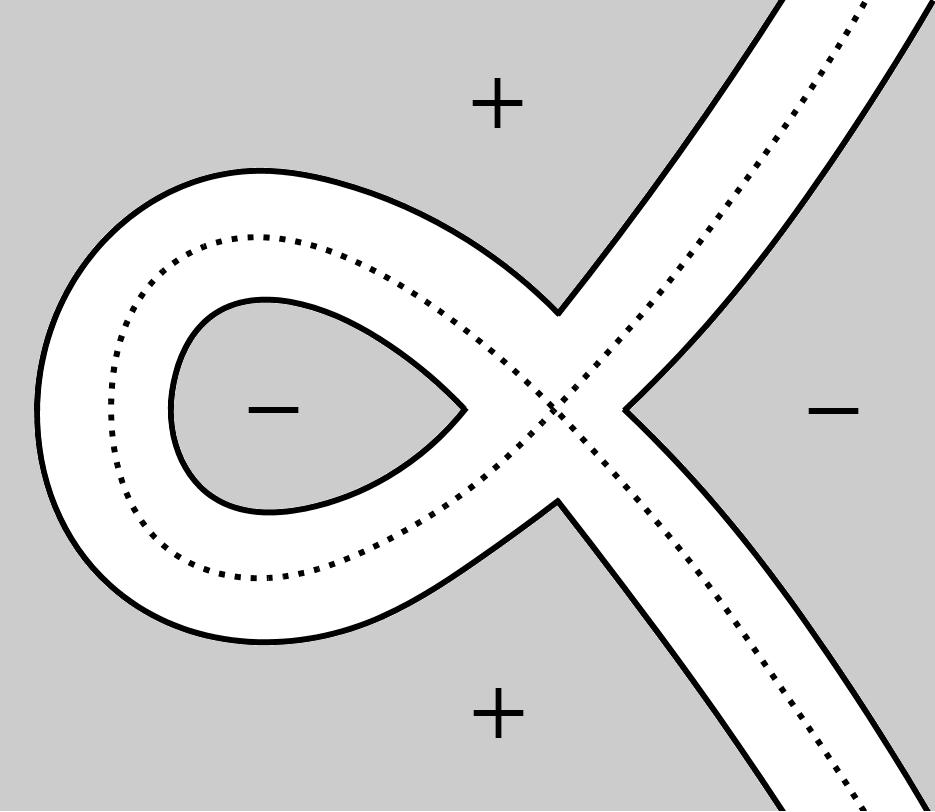}
		\end{minipage}\hfill
		\begin{minipage}[l]{0.5\textwidth}
			\caption{The real-oriented blow-up $\Bl^\R (X, X_0)$ for $\pi(x,y) = y^2 - x^3 - 3x^2$ is represented by the grey regions. 
			The signs $+$ and $-$ indicate 
			         the regions that map to $[0,1)$ and $(-1,0]$ under $\Blpi$, respectively.}
			\label{blowup}
		\end{minipage}		
	\end{figure}%
\end{example}

Our main result in the framework of smooth semi-stable degenerations is the following theorem. We use the shorthand
$0^+ = (+, 0) \in \BlI$ for the positive copy of zero in $\BlI$. The positive half-interval in $\BlI$ 
is consequently denoted by $[0^+,1) \subset \BlI$ for emphasis. 

\begin{theorem} \label{thm:MainSmooth}
  Let $\pi \colon X \to I$ be a smooth semi-stable degeneration and let  
	$\Blpi \colon \BlX \to \BlI$ be the real-oriented blow-up of $\pi$. 
	We set $\BlX^+ := \Blpi^{-1}([0^+,1))$ and $\BlX_0^+ := \Blpi^{-1}(0^+)$. Then there exists a homeomorphism
	\begin{equation} 
		H \colon \BlX_0^+ \times [0,1) \to \BlX^+
	\end{equation}
	such that $\Blpi(H(x,t)) = t$ for all $(x,t) \in \BlX_0^+ \times [0,1)$. 	
	In particular, $\BlX_0^+$ and $X_t$ are homeomorphic for all $t \in (0,1)$.
\end{theorem}

\section{Straightening corners} \label{straightening}

In order to prove \autoref{thm:MainSmooth}, it is convenient to use the language of manifolds
with corners and the well-known (though not well-documented) method of straightening corners.

The first systematic study of manifolds with corners
appears in \cite{Cer-TopologieDeCertains, Dou-VarietesABord}.
We will also use the references \cite{Mel-DifferentialAnalysisManifolds, Joy-ManifoldsCorners}.
A \emph{(smooth) manifold with corners} is a paracompact Hausdorff space $X$ covered by 
a maximal atlas of charts which are homeomorphisms between open sets in $X$ and open sets in
$\Rpos^k \times \R^{n-k}$ containing $0$. We call $(k,n)$ the type of the chart.
The transition maps between charts are required to be smooth. Here, a smooth function
on an open set $U \subset \Rpos^k \times \R^{n-k}$ is the restriction of a smooth function 
on an open set in $\R^n$ containing $U$. 
In particular, this guarantees that each point $p \in X$ has a well-defined
tangent space $T_p X$ isomorphic to $\R^n$. 
For more details, we refer to \cite[Definition 1.6.1]{Mel-DifferentialAnalysisManifolds}
and \cite[Section 2]{Joy-ManifoldsCorners}. Note that \cite{Mel-DifferentialAnalysisManifolds}
calls the  object of our definition \emph{t-manifolds} and requires an additional global condition 
in his notion of manifold with corners.

A \emph{(strictly) inward-pointing vector field} on a manifold with corners $X$ is a smooth vector field $V$
which, with respect to any chart of type $(k,n)$, satisfies $V(0)_i \geq 0$ (or $V(0)_i > 0$, respectively) for all $i = 1, \dots, k$
(cf.\ \cite[Definition 1.13.10]{Mel-DifferentialAnalysisManifolds}, \cite[Page 4]{Joy-ManifoldsCorners}). 
By \cite[Corollary 1.13.1]{Mel-DifferentialAnalysisManifolds}, an inward-pointing vector field can be integrated on compacts, that is, 
for any compact $K \subset X$ there exists $\varepsilon > 0$ and a smooth map
\[
  F \colon K \times [0,\varepsilon) \to X
\]
such that $F(p,0) = p$ and $d F (d/dt) = V$.

A \emph{total boundary defining function} on $X$ is a smooth function $\rho \colon X \to \Rpos$ 
such that $\partial X = \rho^{-1}(0)$ and for each point $p \in \partial X$ there exists a chart of type $(k,n)$ centred in $p$ 
such that $\rho(x) = x_1 \cdots x_k$ in local coordinates
(cf.\ \cite[Lemma 1.6.2]{Mel-DifferentialAnalysisManifolds}, \cite[Definition 2.14]{Joy-ManifoldsCorners}).
Combining the flow of a strictly inward-pointing vector field with the level sets of a total boundary defining function,
we obtain a way of smoothing the corners and turning  $X$ into manifold with boundary.
This construction can be found in various flavours and levels of details in the literature, for example, in 
\cite[Section 2.6, in particular 2.6.4]{Wal-DifferentialTopology} for $n=2$ (where additionally the uniqueness of the construction is discussed). 
For our purposes, the following statement is sufficient. 

\begin{theorem} \label{thm:StraighteningCorners}
  Let $X$ be a manifold with corners and let $\rho \colon X \to \Rpos$ be a proper total boundary defining function. 
	Then there exist $\varepsilon > 0$ and a 
	(continuous) embedding
	\begin{equation} 
		H \colon \partial X \times [0,\varepsilon] \to X
	\end{equation}
	such that $\text{pr}_2 = \rho \circ H$. 
	In particular, $\partial X$ and $\rho^{-1}(t)$ are homeomorphic for all $t \in (0,\varepsilon]$.
\end{theorem}

\begin{proof}
  Since $\rho$ is proper, $\partial X = \rho^{-1}(0)$ is compact. 
	Fix a strictly inward-pointing vector field $V$ on $X$.
	For existence, we can use to the existence 
	of total boundary defining functions, see 
	\cite[Lemma 1.6.2]{Mel-DifferentialAnalysisManifolds}, whose gradients
	are strictly inward-pointing. More explicitly, we can argue as follows:
	First note that strictly inward-pointing vector fields exist locally in a chart. 
	Moreover, sums and positive multiples of inward-pointing vector fields 
	are inward-pointing. Hence using a partition of unity on $X$ 
	(see \cite[Section 1.3 and Lemma 1.6.1]{Mel-DifferentialAnalysisManifolds}), 
	we can construct an inward-pointing vector field on $X$. 
	
	We denote by $F \colon \partial X \times [0,\varepsilon_1) \to X$ 
	the integral flow associated to $V$. 
	Since $V$ is non-zero on $\partial X$, the image $F(\partial X \times [0,\varepsilon_1))$ is a neighbourhood of $\partial X$. 
	Since $\rho$ is proper, there exists $\varepsilon_2 > 0$ such that $\rho^{-1}([0, \varepsilon_2])$ is contained in this neighbourhood. 
	Since $V$ is strictly inward-pointing and $\rho$ is a total boundary defining function, $\rho$ restricted to a flow line starting at $p \in \partial X$
	is strictly increasing locally at $t=0$. Hence, there exists $0 < \varepsilon \leq \varepsilon_2$ such that
	$\rho$ restricted to a flow line is strictly increasing from value $0$ to value $\varepsilon$. 
	We set $K = \rho^{-1}([0, \varepsilon])$ and denote by $f \colon K \to \partial X$ the map which sends a point to the initial point of its flow line. 
	We set $G = (f, \rho) \colon K \to \partial X \times [0, \varepsilon]$. Clearly, $G$ is continuous. Moreover, it is bijective by our assumption
	that $\rho$ is strictly increasing on flow lines up to value $\epsilon$. Since $K$ is compact and $\partial X \times [0, \varepsilon]$ is Hausdorff, 
	we conclude that $H = G^{-1}$ is also continuous. Hence $H \colon \partial X \times [0, \varepsilon] \to K \subset X$ is the embedding required in the statement. 
\end{proof}

\begin{remark} 
  Let us comment on why we refer to the above embedding $H$ as \enquote{straightening corners}. 
	Fix a smooth strictly increasing function $s \colon \Rpos \to \Rpos$ such that 
	$s(0) = \epsilon/2$ and $s(t) = t$ for $t \geq \epsilon$. 
	Using $H$, we can construct a homeomorphism $S \colon X \to \rho^{-1}([\epsilon/2, \infty))$
	such that $\rho \circ S = s \circ \rho$. Via $S$ we can equip $X$ with the
	structure of a manifold with boundary which coincides with the original structure 
	away from the corner locus ($k \geq 2$).  
\end{remark}

We now explain how to apply straightening of corners to prove \autoref{thm:MainSmooth}.
A map $f \colon X \to Y$ between manifolds with corners is called \emph{smooth}
if pull-backs of smooth functions are smooth. 
This is the naive definition, see e.g.\ \cite[Equation (1.10.18)]{Mel-DifferentialAnalysisManifolds}.
In the literature various refined notions exist, e.g.\ \cite[Section 1.12]{Mel-DifferentialAnalysisManifolds}, \cite[Section 3]{Joy-ManifoldsCorners}. 
Given a smooth map $f \colon X \to Y$, we have induced differential maps $df_p \colon
T_p X \to T_q Y$ for all $p \in X$ defined as usual, and $f$ is called an \emph{immersion}
if $df_p$ is injective for all $p \in X$. 

\begin{lemma} \label{lem:CornerStructure}
  Let $X$ be a smooth manifold and $D \subset X$ a normal crossing divisor. Then the real-oriented blow-up
	$\Bl^\R (X,D)$ carries a unique structure of manifold with corners such that
	$\alpha \colon \Bl^\R (X,D) \to X$ is an immersion. 
\end{lemma}

\begin{proof}
  Clearly, there exists a unique smooth structure on $\alpha^{-1}(X \setminus D)$ which satisfies the condition
	of the statement, namely the one induced from $X \setminus D$ via $\alpha$. 
	Moreover, the standard smooth structure on $\Rpos^k \times \R^{n-k}$ is  
	the only one such that 
	\[
	   \Rpos^k \times \R^{n-k} \hookrightarrow \R^k \times \R^{n-k}
	\]
	is an immersion. By the description in \autoref{lem:localblowup} we can hence cover  $\Bl^\R (X,D)$
	by open sets on which existence and uniqueness of such a structure holds. 
	Finally, it is shown in \cite[Proposition 5.2]{HPV-CompactificationHenonMappings} that the coordinate changes are smooth maps 
	(between open sets in $\Rpos^k \times \R^{n-k}$). Hence the result follows. 
\end{proof}

\begin{lemma} \label{lem:BoundaryFunction}
  Let $\pi \colon X \to I$ be a smooth semi-stable degeneration. Then the real-oriented blow-up
	$\Blpi \colon \BlX \to \BlI$ 
	is a smooth map between manifolds with corners. 
	Moreover, 
	the restriction 
	\begin{equation} 
	  \Blpi^+ \colon \BlX^{+} = \Blpi^{-1}([0^+,1)) \to [0,1)		
	\end{equation} 
	is a total boundary defining function for $\BlX^+$
	with boundary $\partial \BlX^+ = \BlX^+_0$.
\end{lemma}

\begin{proof}
  Both statements are immediate consequences of the local description of $\Blpi$ given 
	in \autoref{eq:Blpi}. 
\end{proof}

\begin{example} 
  The real-oriented blow-up $\Bl^\R (X, X_0)$ for $\pi(x,y) = y^2 - x^3 - 3x^2$ from 
	\autoref{ex:blowup} is a two-dimensional manifold with corners with $4$ corner points,
	see \autoref{blowup}.
	In particular, note that our definition allows boundary components to self-intersect as in this example. 
\end{example}

\begin{proof}[\autoref{thm:MainSmooth}]
  By \autoref{lem:CornerStructure}, $\BlX^+$ carries the structure of a manifold with corners 
	with boundary $\partial \BlX^+ = \BlX_0^+$. By \autoref{lem:BoundaryFunction}, $\rho = \Blpi^+ \colon \BlX^+ \to [0,1)$ 
	is a proper total boundary defining function on $\BlX^+$. Let $H_1 \colon \BlX_0^+ \times [0,\varepsilon] \to \BlX^+$ be an embedding as constructed in  \autoref{thm:StraighteningCorners}.
	Note that over $[\epsilon, 1) \subset (0,1)$, $\rho$ is proper and without critical values. Hence $\rho$ can be trivialised by Ehresmann's theorem
	\cite[Theorem 8.5.10]{Dun-ShortCourseDifferential}, 
	that is, there exists a diffeomorphism
	$H_2 \colon X_\epsilon \times [\epsilon, 1) \to \rho^{-1}([\epsilon, 1))$ such that $\text{pr}_2 = \rho \circ H_2$ and $H_2(p,\varepsilon) = p$ for all $p$. 
	Finally $H \colon \BlX_0^+ \times [0,1) \to \BlX^+$ can be constructed using
	\begin{equation} 
		         (p,t) \mapsto 
						\begin{cases}
						  H_1(p,t) & \text{ if } t \leq \varepsilon, \\
							H_2(H_1(p,\varepsilon), t) & \text{ if } t \geq \varepsilon. 
						\end{cases}
	\end{equation}
\end{proof}

\section{Real semi-stable degenerations}

We now return to the complex and real setting presented in the introduction.
A notational warning: In this context, we will denote by $(X,D)$
the \emph{complex} semi-stable degeneration. 
After adding a suitable real structure, the real 
sublocus $(\R X, \R D)$ forms a \emph{smooth} semi-stable degeneration
in the sense of the previous sections. 
So, the tuple $(X,D)$ from the previous sections will be replaced by
$(\R X, \R D)$ now.

\begin{definition} \label{def:complexSemiStable}
  A \emph{complex semi-stable degeneration} is a proper holomorphic function $\pi \colon X \to \DD$ 
	from a complex manifold $X$ to the unit disc $\DD \subset \C$ such that 
	\begin{itemize}
		\item the map $\pi$ is regular at all $t \in \DD^* = \DD \setminus \{0\}$,
		\item the special fibre $X_0 = \phi^{-1}(0)$ is a reduced normal crossing divisor. 
	\end{itemize}
\end{definition}

The notions of normal crossing divisor\fshyp{}function as given in \autoref{SmoothSemiStable}
admit obvious generalizations to the holomorphic setting. It is then clear that 
\autoref{def:complexSemiStable} and the holomorphic version of \autoref{def:smoothSemiStable}
agree: For any point $p \in X_0$, there exists an holomorphic chart
centred in $p$ such that $\pi(z) = z_1 \cdots z_k$
in local coordinates for some $1 \leq k \leq n$. 
We denote the tangent cone of $X_0$ at $p$ by $T_p X_0$. 
With respect to the chart from above, $T_p X_0 = \{z_1 \cdots z_k = 0\} \subset \C^n = T_p X$ is 
just the union of the first $k$ hyperplanes.

\begin{definition} \label{def:realSemiStable}
  A \emph{real structure} for $\pi \colon X \to \DD$ is an orientation-reversing involution $\iota \colon X \to X$
	such that $\text{conj} \circ \pi = \pi \circ \iota$. 
	The tuple $(\pi, \iota)$ is a \emph{real semi-stable degeneration}. We say that $(\pi, \iota)$ 
	is \emph{totally real} if 
	for any $p \in \R X_0$ the differential $d \iota \colon T_p X \to T_p X$ 
	fixes each hyperplane in the tangent cone $T_p X_0$.	
\end{definition}

Throughout the following, the letter $\R$ indicates restriction to real points, that is, to the fixed points 
of a (given) real structure $\iota$.

\begin{lemma} \label{lem:TotallyRealToSmooth}
  If $(\pi, \iota)$ is totally real, then the map $\R\pi \colon \R X \to \R \DD = I = (-1, 1)$ is a smooth semi-stable degeneration. 
\end{lemma}

\begin{proof}
  The definition implies that for any $p \in \R X_0$ there exists a holomorphic normal crossing chart $U \subset \C^n \to \C$, $\pi(z) = z_1 \cdots z_k$
	for $\pi$ at $p$ such that $\iota$ corresponds to the standard (coordinate-wise) conjugation in $\C^n$. Clearly, the real locus of this chart 
	$U \cap \R^n \to \R$, $\R \pi(x) = x_1 \cdots x_k$ can be used as smooth normal crossing chart for 
	$\R\pi \colon \R X \to \R \DD$.
\end{proof}

\begin{example} \label{ex:cubiccomplex}
  The smooth semi-stable degeneration from \autoref{ex:cubic} is the real locus of the totally real semi-stable degeneration
	given by the same map $\pi(x,y) = y^2 - x^3 -3 x^2$ allowing complex arguments and values (in $\DD$). 
\end{example}

\begin{example} \label{ex:smallperturbation}
  \autoref{ex:cubic}\fshyp{}\autoref{ex:cubiccomplex} is a special case of the classical \enquote{small perturbation} method 
	that was used in particular by Harnack and Hilbert
	to construct real projective curves with interesting topological properties, see
	\cite{Har-UeberDieVieltheiligkeit,Hil-UeberDieReellen} and \autoref{smallperturbation}.
	We refer to \cite[Theorem 1.5.A]{Vir-IntroductionTopologyReal}
	for more background. The construction can in fact be regarded as an example	of semi-stable degenerations. 
	The input for the small perturbation consists of 
	two real projective curves $C$ and $D$ of degree $d$ such that
  \begin{itemize}
		\item the singularities of the real part $\R C$ are non-isolated nodes (that is, locally given by $x^2-y^2 =0$),
		\item the real parts $\R C$ and $\R D$ intersect transversally (in particular, in non-singular points of $C$ and $D$).
	\end{itemize}
	We now consider the pencil of curves $C_t = V(F + t G)$ 
	where $F$ and $G$ are homogeneous polynomials 
	defining $C$ and $D$, respectively.
	The method is then based on describing $\R C_t$ for small $t$
	in terms of $\R C$ and certain combinatorial data (the pattern of signs
	of $FG$ on $\R\P^2 \setminus (\R C \cup \R D)$). Note that the conditions
	on the input data are essentially equivalent to the (smooth) semi-stability of the 
	family 
	\[
	  \{F + tG = 0\} \subset \R\P^2 \times I \to I.
	\]
	Moreover, in practice one can often assume that 
	the conditions extend to complex points, that is, that additionally $C$ is nodal 
	and $C$ and $D$ intersect transversally. In this case, the above 
	degeneration is the real locus of the complex semi-stable degeneration
	\[
		\{F + tG = 0\} \subset \C\P^2 \times \DD \to \DD.
	\]
	The requirement for the nodes of $\R C$ to be non-isolated 
	is equivalent to the degeneration being totally real. 

  \begin{figure}[t]
	\begin{minipage}[r]{0.4\textwidth}
		\includegraphics[width=\textwidth]{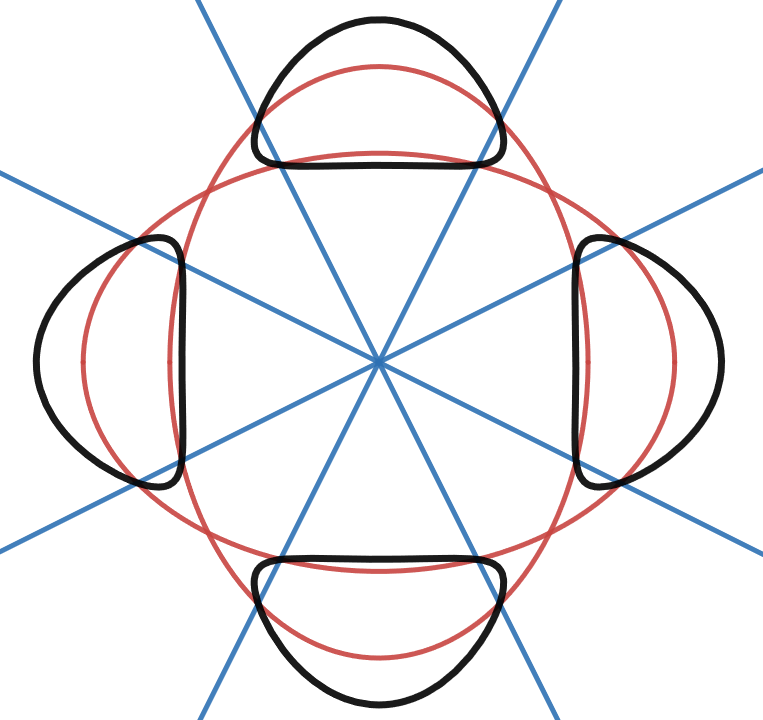}
	\end{minipage}\hfill
	\begin{minipage}[l]{0.5\textwidth}
		\caption{The small perturbation $C_t$ of $C$, the union of two conics, along $D$, the union of four lines. The result is a curve of degree $4$
		with $4$ connected components.}
		\label{smallperturbation}
	\end{minipage}		
  \end{figure}%
	
	The present work can hence be regarded as an attempt to generalize the small perturbation method
	to arbitrary dimensions and the abstract setting. 
	Note, however, that our methods would have to be modified in order to describe
	the topological pair $\R C_t \subset \R\P^2$ (not just $\R C_t$). One approach would be to 
	also degenerate $\R\P^2$ and use the relative versions of our statements, c.f.\ \autoref{toricdegenerations}. 
\end{example}

\begin{example} \label{ex:normalcone}
  A classical example of semi-stable degenerations are deformations to the normal cone \cite[Chapter 5]{Ful-IntersectionTheory}. 
	Let $Y$ by a compact complex manifold with real structure
	and $Z \subset Y$ a real submanifold. We denote by $X$ the (ordinary) blow-up of $Y \times \DD$ along $Z \times \{0\}$. Then the canonical projection map
	$\pi \colon X \to \DD$ is a totally real semi-stable degeneration. The generic fibre is $Y$, while the special fibre consists of $B = \Bl_Z(Y)$ and
	the projective completion $C$ of the normal bundle $N_Z Y$ of $Z$ in $Y$, that is, $C = \P(N_Z Y \oplus \OO_Z)$. These two components intersect 
	in the exceptional divisor of $B$ and the section at infinity of $C$ which is the projective normal bundle $E = \P(N_Z Y)$.
	Note that if $Z \subset Y$ is a divisor, then $B = Y$ and $E = Z$.
\end{example}

Given a normal crossing divisor $D \subset X$ in a complex manifold, we can apply 
the more general definition of real-oriented blow-up $\Bl^\R(X,D)$ given 
in \cite{HPV-CompactificationHenonMappings}  to $D \subset X$ regarded
as real-analytic (sub-)varieties. However, the outcome has 
a simple description in terms of polar coordinates 
which allows us to avoid discussing the general definition.
Again, it is sufficient to understand the local case. 
We denote by $S^1$ the unit circle in $\C$. 

\begin{lemma} \label{lem:localcomplexblowup}
  Set $X = \C^n$ and $D = \{z_1 \cdots z_k = 0\}$ for some $1 \leq k \leq n$. Then
	the real-oriented blow-up of $X$ along $D$ is given by
	\begin{equation} 
	  (S^1 \times \Rpos)^k \times \C^{n-k} \to \C^k \times \C^{n-k}
	\end{equation}
	where the map on the first $k$ factors is given by 
	\begin{equation} \begin{split} 
		S^1 \times \Rpos &\to \C, \\
		(\theta, r) &\mapsto \theta r	
	\end{split} \end{equation}
	and on the last $n-k$ factors is identity. 
\end{lemma}

\begin{proof}
  This follows from \cite[Theorem 5.3]{HPV-CompactificationHenonMappings} and the base case
	\[
	  \Bl^\R (\C, \{0\}) = S^1 \times \Rpos \to \C
	\]
	discussed for example at \cite[Page 235]{HPV-CompactificationHenonMappings}.  
\end{proof}

Again, the construction of real-oriented blow-ups can be applied 
to (suitable) families of complex varieties. 
This is discussed in detail in 
\cite[Chapter X, §9 and Chapter XV, §8]{ACG-GeometryAlgebraicCurves}
in the case of curves.
More generally, given a complex semi-stable degeneration 
$\pi \colon X \to \DD$, we can construct the real-oriented blow-up
$\Blpi \colon \BlX \to \BlDD$ which locally looks like 
\begin{equation} \begin{split} 
  \Blpi \colon (S^1 \times \Rpos)^k \times \C^{n-k} & \to S^1 \times \Rpos, \\
	              ((\theta, r), z) & \mapsto (\theta_1 \cdots \theta_k, r_1 \cdots, r_k).
\end{split} \end{equation}

\begin{proof}[\autoref{thm:MainComplex}]
	When $\pi \colon X \to \DD$ is totally real, 
	the blown-up family $\Blpi \colon \BlX \to \BlDD$ carries an induced real 
	structure given, in local charts and coordinate-wise, by $S^1 \times \Rpos \to S^1 \times \Rpos$, 
	$(\theta, r) \mapsto (1/\theta, r)$. 
	The relation to the blow-up of
	$\R\pi \colon \R X \to (-1, 1)$ 
	is the expected one: 
	We have $\R\BlX = \widetilde{\R X}$ and $\R\Blpi = \widetilde{\R\pi}$. 
  Hence, \autoref{thm:MainComplex} is the consequence of applying \autoref{thm:MainSmooth}
	to the smooth semi-stable degeneration $\R\pi \colon \R X \to (-1, 1)$. 
\end{proof}

\section{The relative version}

Given a smooth semi-stable degeneration $\pi \colon X \to I$, 
for any point $p \in X$ with $t = \pi(p)$ the tangent cone $T_p X_t$ is a union 
of a finite number of hyperplanes $H_{p,1}, \dots, H_{p,k} \subset T_p X$. 
Note that $k = 1$ if $t \neq 0$ or if $p$ is a generic point of $X_0$. 
A collection of closed submanifolds $X_1, \dots, X_m \subset X$
is called \emph{transversal} if for each $p \in X$ 
the collection of subspaces $H_{p,1}, \dots, H_{p,k}, T_p X_1, \dots, T_p X_m \subset T_p X$
is transversal. That is, we require the codimension of 
\[
  H_{p,1} \cap \dots \cap H_{p,k} \cap T_p X_1 \cap \dots \cap T_p X_m \subset T_p X
\]
to be the sum of the codimensions of the individual spaces. 

\begin{theorem} \label{thm:RelativeSmooth}
  Let $\pi \colon X \to I$ be a smooth semi-stable degeneration and 
	$X_1, \dots, X_m \subset X$ a transversal collection of submanifolds. 
	\begin{enumerate}
		\item For any $I \subset \{1, \dots, m\}$ such that the intersection $X_I = \bigcap_{i \in I} X_i$ is non-empty, 
		      the restriction $\pi_I = \pi|_{X_I} \colon X_I \to I$ is a semi-stable degeneration. 
		\item There exists a homeomorphism $H \colon \BlX_0^+ \times [0,1) \to \BlX^+$ as in \autoref{thm:MainSmooth}
		      such that for each $I \subset \{1, \dots, m\}$ we have 
					\[
					  H\left(\BlX^+_{I,0} \times [0,1)\right) = \BlX^+_{I}. 
					\]
	\end{enumerate}
\end{theorem}

\begin{proof}
  Part (a): Assume that $X_I$ is non-empty. Since $X_I$ is closed in $X$, the restriction $\pi_I$ is a proper map. 
	          Transversality of the $T_p X_1, \dots, T_p X_m$ ensures that $X_I$ is a submanifold.
						Transversality of $T_p X_I$ and $T_p X_t$ for $\pi(p) = t \neq 0$ ensures that $\pi_I$ is regular over $I \setminus \{0\}$. 
						Finally, for $p \in X_I \cap X_0$, we can 
						find numbers $1 \leq k_0 \leq k_1 \leq \dots \leq k_m \leq n$
						and a normal crossing chart for $\pi$ centred in $p$ such that in local coordinates
						\begin{equation} \label{eq:normalchart} 
						  \pi(x) = x_1 \cdots x_{k_0} \text{ and } X_{i} = \{x_{k_{i-1} + 1} = \dots = x_{k_i} = 0\} \text{ for } i = 1, \dots, m.
						\end{equation}
						Here, it is understood that $k_{i-1} = k_i$ means that $p \notin X_i$.
						Hence the restriction of this chart to $X_I$ corresponds to setting some coordinates 
						$x_j$, $j > k_0$, to zero and hence provides a normal crossing chart for $\pi_I$ at $p$. 
						
	Part (b): Using again the normal chart from \autoref{eq:normalchart}, we first observe 
	          that the real-oriented blow-up $\BlX_I$ of $X_I$ is in fact equal to the total transform
						$\alpha^{-1}(X_I) \subset \BlX$. It now suffices to choose the vector field $V$ that is used in the proof 
						of \autoref{thm:StraighteningCorners} in a more specific manner. More precisely, we
						require for any $p \in \BlX$ with $I(p) := \{ i : p \in \BlX_i\}$ that $V(p)$ is contained in $T_p \BlX_{I(p)}$. 
						From \autoref{eq:normalchart}, it is clear that such a vector field exists locally.
						Since the required property is invariant under rescaling, we can use partition of unity
						to obtain a global vector field with the required property. 
						We refer to \cite[Section 1.3]{Mel-DifferentialAnalysisManifolds} and \cite[Section 1.2.6]{Cer-TopologieDeCertains}
						for a discussion of partitions of unity on manifolds with corners. 
						Now, following the procedure in the proof of \autoref{thm:StraighteningCorners} we obtain 
						flow lines for which $I(p)$ is an invariant and hence the map $H_1 \colon \BlX^+_0 \times [0,\epsilon] \to \BlX^+$
						satisfies $H_1(\BlX^+_{I,0} \times [0,\epsilon]) \subset \BlX^+_I$ for all $I$. 
						Combining this with the relative version of Ehresmann's theorem, 
						we can conclude as 
						in the proof of \autoref{thm:MainSmooth}.
\end{proof}

Let $\pi \colon X \to \DD$ be a complex semi-stable degeneration. A collection of 
closed complex submanifolds $X_1, \dots, X_m \subset X$ is called \emph{transversal} 
if it satisfies the same condition as in the smooth case, applied to the complex subspaces
of $T_p X$. 
As in the proof of Part (a) of \autoref{thm:RelativeSmooth}, we conclude that the 
restrictions $\pi_I = \pi|_{X_I} \colon X_I \to I$ are complex semi-stable degenerations. 
With these definitions, \autoref{thm:RelativeReal} is now a simple corollary of 
\autoref{thm:RelativeSmooth}.

\begin{proof}[\autoref{thm:RelativeReal}]
	Given a point $p \in \R X_0$, we may use 
	the complex version of the normal chart in \autoref{eq:normalchart} 
	(with the canonical real structure) to 
	conclude that $\pi_I$ is totally real for all $I \subset \{1, \dots, m\}$
	and that the collection $\R X_1, \dots, \R X_m \subset \R X$ is a transversal collection 
	in the smooth sense. 
	Hence the statement follows from \autoref{thm:RelativeSmooth}.
\end{proof}

\section{Stratifying the special fibre} \label{secStratifying}

Let $X$ be a smooth manifold and $D \subset X$ a normal crossing divisor. 
For a point $p \in X$, any normal crossing chart centred in $p$
is of the same type $(k,n)$. We call $k$ the \emph{codimension} of $p$ in $(X,D)$. 
The \emph{codimension $k$ skeleton} $X^k$ is the set of points 
of codimension at least $k$. In particular $X^0 = X$ and $X^1 = D$.
The successive complements $X^k_\circ = X^k \setminus X^{k+1}$
(the locus of points of codimension equal to $k$) 
are locally closed submanifolds of $X$. 
A connected component $S \subset X^k_\circ$ is called a \emph{stratum}
of codimension $k$ of $(X,D)$. Clearly $X$ is the disjoint union of all these strata. 
Moreover, the strata satisfy the \emph{Whitney condition}, see 
\cite[Part I, Section 1.2]{GM-StratifiedMorseTheory}. Indeed, since the condition is local,
it suffices to check it for a normal crossing chart where it is obvious. 
Hence $X^n \subset \dots X^1 \subset X^0$ defines a \emph{Whitney stratification}
on $X$ in the sense of \cite[Part I, Section 1.2]{GM-StratifiedMorseTheory}.
We denote the set of strata by $\SS(X,D)$. 

From now on, we denote by $\SS$ a Whitney stratification that refines $\SS(X,D)$. 
In  our applications, the refinement $\SS$ will be equal to 
$\SS(X,D')$ where $D'$ is a normal crossing divisor
containing $D$. 
For any $S \in \SS$, we denote by $k(S)$ the codimension of the unique stratum in $\SS(X,D)$
that contains $S$. In other words, $k(S)$ is the number of sheets of $D$ that intersect in 
a point $p \in S$. 
Let $\alpha \colon \BlX \to X$ be the real-oriented blow-up of $X$ along $D$. 
We set $\wt{S} = \alpha^{-1}(S)$ for any $S \in \SS$ and 
denote by $\chi^c$  
the Euler characteristic with closed support. 

\begin{theorem} \label{thm:StrataSmooth}
  The real-oriented blow-up $\BlX = \Bl^\R(X,D)$ of $X$ along a normal crossing divisor $D$ 
	is a Whitney-stratified space via (the connected components of)
	\[
	  \BlX = \bigsqcup_{S \in \SS} \wt{S}.
	\]
	Moreover, 
	\begin{align} 
		\chi^c(\BlX) = \sum_{k=0}^n \sum_{\substack{S \in \SS \\ \cd(S) = k}} 2^k \chi^c(S), & & \chi^c(\wt{D}) = \sum_{k=1}^n \sum_{\substack{S \in \SS \\ \cd(S) = k}} 2^k \chi^c(S).
	\end{align}
\end{theorem}

\begin{proof}
	Any manifold with corners, hence also $\BlX$, carries a canonical 
	Whitney stratification (given again by codimension of points) denoted 
	by $\SS(\BlX)$. 
	The map $\alpha$ is a stratified map with respect to $\SS(X,D)$ and 
	$\SS(\BlX)$, that is, the image of a stratum in $\SS(\BlX)$
	is a stratum in $\SS(X,D)$ and the restriction of $\alpha$ to these strata
	is a submersion. In fact, in our situation such a restriction 
	is also an immersion and hence a covering map.
	In particular, the (connected components) of the $\wt{S}$ form a refinement 
	of $\SS(\BlX)$. 
	Moreover, the restriction $\alpha_S \colon \wt{S} \to S$ is a smooth covering map
	of degree $2^k$ with $k = \cd(S)$. 
	All these statements can be checked immediately using the local description of 
	$\alpha$ in \autoref{eq:localblowup}. Hence
	\begin{align*} 
		\chi^c(\BlX) = \sum_{S \in \SS}  \chi^c(\wt{S}) & & \text{and} & & \chi^c(\wt{S}) = 2^k \chi^c(S)
	\end{align*}
	for $S \in \SS$ with $k = \cd(S)$, which proves the theorem. 
\end{proof}

\begin{example} 
  \autoref{strata} depicts the stratification of $D$ and $\wt{D}$ from \autoref{ex:cubic}.
	\begin{figure}[h]
		\begin{minipage}[r]{0.4\textwidth}
		  \centering
			\includegraphics[width=0.8\textwidth]{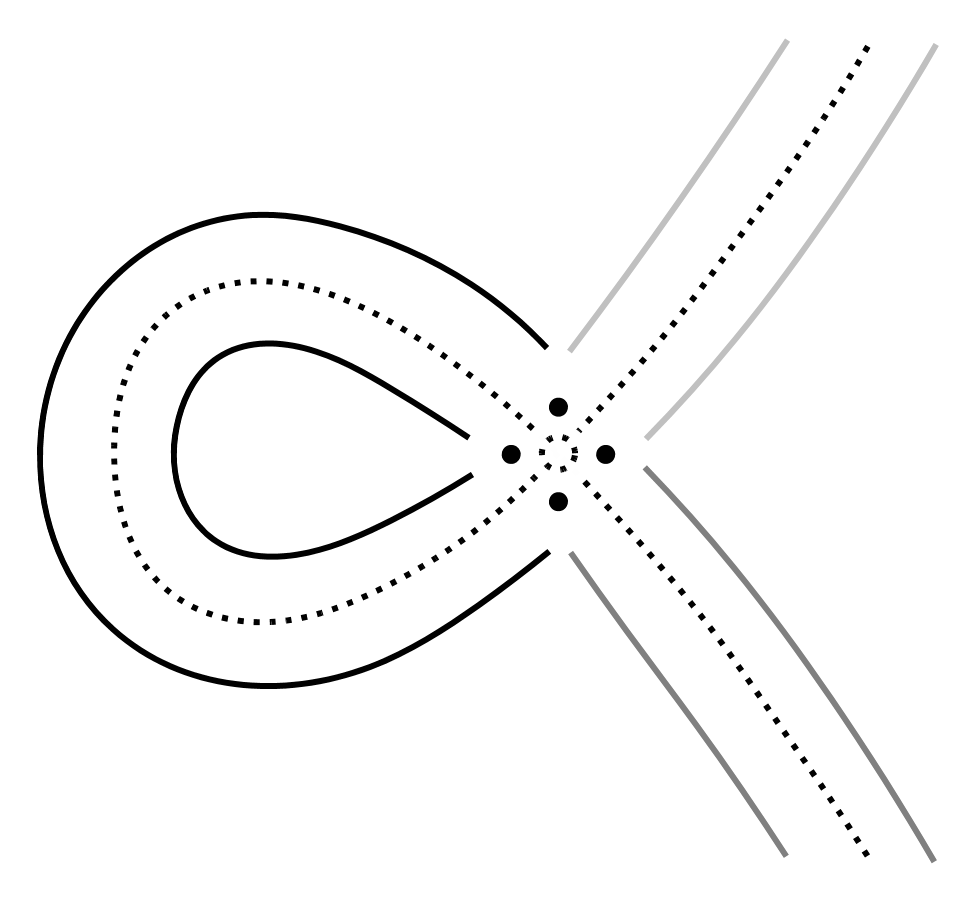}
		\end{minipage}\hfill
		\begin{minipage}[l]{0.5\textwidth}
			\caption{The dashed segments and point depict the stratification of the normal crossing divisor $D$. The solid segments and points
			        represent the strata of $\wt{D}$. All four points belong to $\wt{S}$ where $S= \{(0,0)\}$ is the node. 
							A one-dimensional $\wt{S}$ consists of two segments of the same colour.}
			\label{strata}
		\end{minipage}		
	\end{figure}
\end{example}

In the case of a smooth semi-stable degeneration 
$\pi \colon X \to I$, we can restrict the previous discussion to the 
positive part. For any $S \in \SS$ we set 
$\wt{S}^+ = \wt{S} \cap \BlX^+$.

\begin{corollary} \label{cor:StrataSmoothSemiStable}
	The positive special fibre $\BlX^+_0$ of a smooth semi-stable degeneration 
	$\pi \colon X \to I$ is a Whitney-stratified space via 
	(the connected components of)
	\[
	  \BlX^+_0 = \bigsqcup_{\substack{S \in \SS \\ S \subset X_0}} \wt{S}^+.
	\]
	Moreover, 
	\begin{equation} 
		\chi^c(\BlX^+_0) = \sum_{k=1}^n \sum_{\substack{S \in \SS \\ \cd(S) = k}} 2^{k-1} \chi^c(S).
	\end{equation}
\end{corollary}

\begin{proof}
  The statement follows from \autoref{thm:StrataSmooth} plus the observation
	(see \autoref{eq:Blpi}) 
	that for any $p \in X_0$, the subset of preimages $q \in \alpha^{-1}(p)$ such that 
	$\Blpi(q) = 0^+$ consists of exactly half the points, that is, $2^{k-1}$ elements.   
\end{proof}

\begin{proof}[\autoref{thm:StrataComplex}]
  The statement is a consequence of applying \autoref{cor:StrataSmoothSemiStable} 
	to the smooth semi-stable degeneration $\R \pi \colon \R X \to \R \DD$ 
	from \autoref{lem:TotallyRealToSmooth}.
\end{proof}

\begin{example} \label{ex:doublecovering}
  Using the notation from \autoref{ex:normalcone}, 
	let us consider the semi-stable degeneration given 
	by the deformation of $Y = \C\P^2$ to the normal cone 
	of a line $Z \subset \C\P^2$.
	In this case the special fibre $X_0$ consists of the two components
	$B = \C\P^2$ and $C = \P(\OO_Z(1) \oplus \OO_Z) = \Bl_p \C\P^2$ glued
	along $E = Z = \C\P^1$. 
	Hence the stratification of $\R X_0$ consists of the three strata
	$S_E = \R E = S^1$, $S_B = \R B \setminus \R E = \R^2$ and $S_C = \R C \setminus \R E = \R \P^2 \setminus \{\text{pt}\}$.
	Clearly, we have $\wt{S}^+_B = S_B$ and $\wt{S}^+_C = S_C$ since $k=1$. 
	The map $\wt{S}^+_E \to S_E$, however, is the \emph{connected} degree $2$ covering map of $S^1$ (here $k=2$). 
	Indeed, note that the normal bundles $N_{\R E} \R B$ and $N_{\R E} \R C$ are Möbius strips. 
	The Euler characteristic computation in this case reads as 
	\begin{align*} 
		\chi( \R \P^2) & = \chi^c(S_B) + \chi^c(S_C) + 2 \chi(S_E) \\
		               & = \chi^c(\R^2) + \chi^c(\R \P^2 \setminus \{\text{pt}\}) + 2 \chi(S^1)
									 & = 1 + 0 + 0 = 1.
	\end{align*}
\end{example}

\begin{remark} \label{rem:weaklySemiStable}
  At this point, it might be worthwhile to mention that most of what we did
	also works under the weaker assumption that $\pi$ is a weakly semi-stable.
	Here, a \emph{weakly semi-stable degeneration} $\pi$ is a degeneration 
	such that the total space $X$ is smooth, $\pi$ is regular away from $0$, and 
	$X_0 \subset X$ is a normal crossing divisor, and $\pi$ is of finite order at $X_0$
	(e.g.\ $\pi$ is real-analytic), but not necessarily normal crossing, that is, not necessarily reduced at $X_0$. 
	In other words, a weakly semi-stable degeneration $\pi$ has the local form
	\[
	  \pi(x_1, \dots, x_n) = \pm x_1^{a_1} \dots x_k^{a_k}
	\]
	for some positive integers $a_1, \dots, a_k \geq 1$. 
	Using these local charts, it is straightforward to check that all statements
	made until (including) \autoref{thm:StrataSmooth} remain true without change.
	\autoref{cor:StrataSmoothSemiStable} can be adapted to weakly semi-stable degenerations as follows:
	The local form from above can be subdivided into three cases,
	\begin{align} 
		 \pi(x_1, \dots, x_n) &= x_1^{a_1} \dots x_k^{a_k}, & & \text{at least one $a_i$ is odd}, \nonumber \\
		 \pi(x_1, \dots, x_n) &= x_1^{a_1} \dots x_k^{a_k}, & & \text{all $a_i$ are even}, \label{eq:weaklyTypes} \\
		 \pi(x_1, \dots, x_n) &= -x_1^{a_1} \dots x_k^{a_k}, & & \text{all $a_i$ are even}. \nonumber 
	\end{align}
	Given a stratum $S \subset X_0$, the covering map $\wt{S}^+ \to S$ is 
	of degree $2^{k-1}$, $2^k$ or $0$, respectively, depending on the associated local form
	in the above order. 
\end{remark}

\section{Gluing the special fibre} \label{describing}

Let $X$ be a smooth manifold with normal crossing divisor $D$. 
In applications, it is of interest to provide a more explicit description of 
how the various strata are glued together to form $\BlX$ and $\wt{D}$, respectively. 
To do so, it is useful to introduce a few more definitions. 
Fix $S \in \SS$ with $k = \cd(S)$.  
For any $p \in S$, the preimage $\alpha^{-1}(p)$ can be canonically identified
with the \emph{orthants near $p$} given by 
\begin{equation} 
  Q(p) = \varprojlim_{U \ni p \text{ open}} \pi_0 (U \setminus D) = \pi_0(T_p X \setminus T_p D).
\end{equation}
Given a chart of type $(k,n)$ at $p$, we get a canonical bijection between 
$Q(p)$  and $\pi_0 (\R^n \setminus \{x_1 \cdots x_k = 0\})$ which then 
yields a bijection to $\{\pm\}^k$. 
In particular, $|Q(p)| = 2^k = |\alpha^{-1}(p)|$.
In the case of a smooth semi-stable degeneration $\pi \colon X \to I$,
we can additionally define the subset of \emph{positive} orthants
$Q^+(p) \subset Q(p)$ given by 
\begin{equation} 
  Q^+(p) = \varprojlim_{U \ni p \text{ open}} \pi_0 (U \cap \pi^{-1}((0,1))).
\end{equation}
For $p \in X_0$, this subset can be canonically identified with $\alpha^{-1}(p) \cap \Blpi^{-1}(0^+)$
and hence $|Q^+(p)| = 2^{k-1}$.

In general, in order to describe how the various strata are glued together, 
the monodromy of the covering maps $\wt{S} \to S$ has to be taken into account.
For simplicity, we restrict our discussion here to the case of simply connected strata 
(or if we know for some other reason that the monodromy is trivial). 
A stratification $\SS$ is called \emph{simply connected} if all its strata are simply connected. 
From now we assume that $\SS$ is simply connected. 
This implies, of course, that $\wt{S} = S \times Q(S)$ is a trivial covering map
of $S$ whose sheets are labelled by $Q(S)$. Moreover, in the semi-stable degeneration
case, we have $\wt{S}^+ = S \times Q^+(S)$.

Given a stratum $S \in \SS$ and a point $p \in \overline{S}$, we define the \emph{ends of $S$ in $p$} by
\begin{equation} 
  E(p,S) = \varprojlim_{U \ni p \text{ open}} \pi_0 (U \cap S).
\end{equation}
Given another point $p'$ in the same stratum $T \in \SS$ as $p$ and a continuous path connecting
the two, we get induced identifications of $E(p,S)$ and $E(p',S)$ as well as 
$Q(p)$ and $Q(p')$. Since $T$ is simply connected, these identifications do no depend on the path. 
Hence we allow ourselves to write $E(T,S)$, $Q(S)$ and $Q^+(S)$ in the following. 
Note that $|E(T,S)| \neq 1$ corresponds to self-intersections of $\overline{S}$ in $T$, 
c.f.\ \autoref{blowup} and \autoref{gluing}.
We have a canonical \emph{orthant map}
\begin{equation} 
  \mu_{T,S} \colon E(T,S) \times Q(S) \to Q(T)
\end{equation} 
defined as follows: Consider a normal crossing chart
of type $(k,n)$ centred at a point in $T$. Via such a chart, $Q(T)$ corresponds
to the orthants of $\R^k$. An end $e \in E(T,S)$ corresponds to a face $F$
of the subdivision of $\R^k$ given by the coordinate hyperplanes. Finally, the orthants
adjacent to $F$ can be canonically identified with $Q(S)$. Then 
$\mu(e, .) \colon Q(S) \hookrightarrow Q(T)$ is the map corresponding 
to the inclusion.

The \emph{endpoint modification} of $S$ 
is the map $\beta_S \colon \widehat{S} \to \overline{S} \subset X$ whose 
fibre over $p \in X$ is $\beta_S^{-1}(p) = E(p,S)$. 
We can equip $\widehat{S}$ with a natural topology by requiring
that $\beta_S$ is continuous and a sequence $p_k \in S$ converges
to $e \in E(p,S)$ if it converges to $p$ (in $\overline{S}$) 
and, for each $p \in U$, the $p_k$ lie in the component
of $U \cap S$ specified by $e$, for large $k$. 
We assume from now on that $\SS = \SS(X,D')$ for some normal crossing divisor $D'$ containing
$D$. In this case, using normal crossing charts 
we check easily that $\widehat{S}$ carries the structure of a manifold with corners 
that turns $\beta_S$ into a stratified smooth map with respect to $\SS(\widehat{S})$. 
For more details, we refer to \cite[Page 235]{HPV-CompactificationHenonMappings},
\cite[Definition 2.6]{Joy-ManifoldsCorners} and \cite[Section 6]{Dou-VarietesABord}.
Note that $\widehat{S} \neq \overline{S}$ is equivalent  to $|E(p,S)| \neq 2$ for some $p$ which 
as mentioned above corresponds to self-intersections, see
\autoref{gluing} for an example.

\begin{theorem} \label{thm:GluingSmooth}
  The real-oriented blow-up $\BlX = \Bl^\R(X,D)$ of $X$ along a normal crossing divisor $D$ 
	is the Whitney-stratified space 
	obtained from gluing
	\[
		\widehat{X} = \bigsqcup_{S \in \SS} \widehat{S} \times Q(S)
	\] 
	via the maps
	\begin{equation} \begin{split} 
		\partial\widehat{S} \times Q(S) = \bigsqcup_{T \subsetneq S} T \times E(T,S) \times Q(S) &\to \bigsqcup_{T \subsetneq S} T \times Q(T) \subset \widehat{X}, \\
		(p, e, q) &\mapsto (p, \mu_{T,S}(e,q)).
	\end{split} \end{equation}
	The positive special fibre $\BlX^+_0$ of a smooth semi-stable degeneration 
	$\pi \colon X \to I$ is the  Whitney-stratified space 
	obtained from gluing
	\begin{equation} \label{eq:startingspaceSmooth} 
	  \widehat{X}^+_0 = \bigsqcup_{\substack{S \in \SS \\ S \subset X_0}} \widehat{S} \times Q^+(S)
	\end{equation}	  
	via the same maps restricted to positive orthants. In particular, in this case $\mu_{T,S}(e,q) \in Q^+(T)$ for all
	$e \in E(T,S)$, $q \in Q^+(S)$. 
\end{theorem}

\begin{proof}
  This is a consequence of the preceding discussion, except for possibly the last point:
	In the case of a semi-stable degeneration, given $T \subsetneq S$ 
	we may use a normal crossing chart
	for $\pi$ (of type $(k,n)$) to define $\mu_{T,S}$. Then the positive orthants for both
	$S$ and $T$ are the orthants on which $x_1 \cdots x_k$ is positive. Hence
	$\mu$ respects positivity. 
\end{proof}

\begin{example} 
  The gluing process for \autoref{ex:cubic} is depicted in \autoref{gluing}.
  \begin{figure}[ht]
	\begin{minipage}[r]{0.6\textwidth}
		\includegraphics[width=\textwidth]{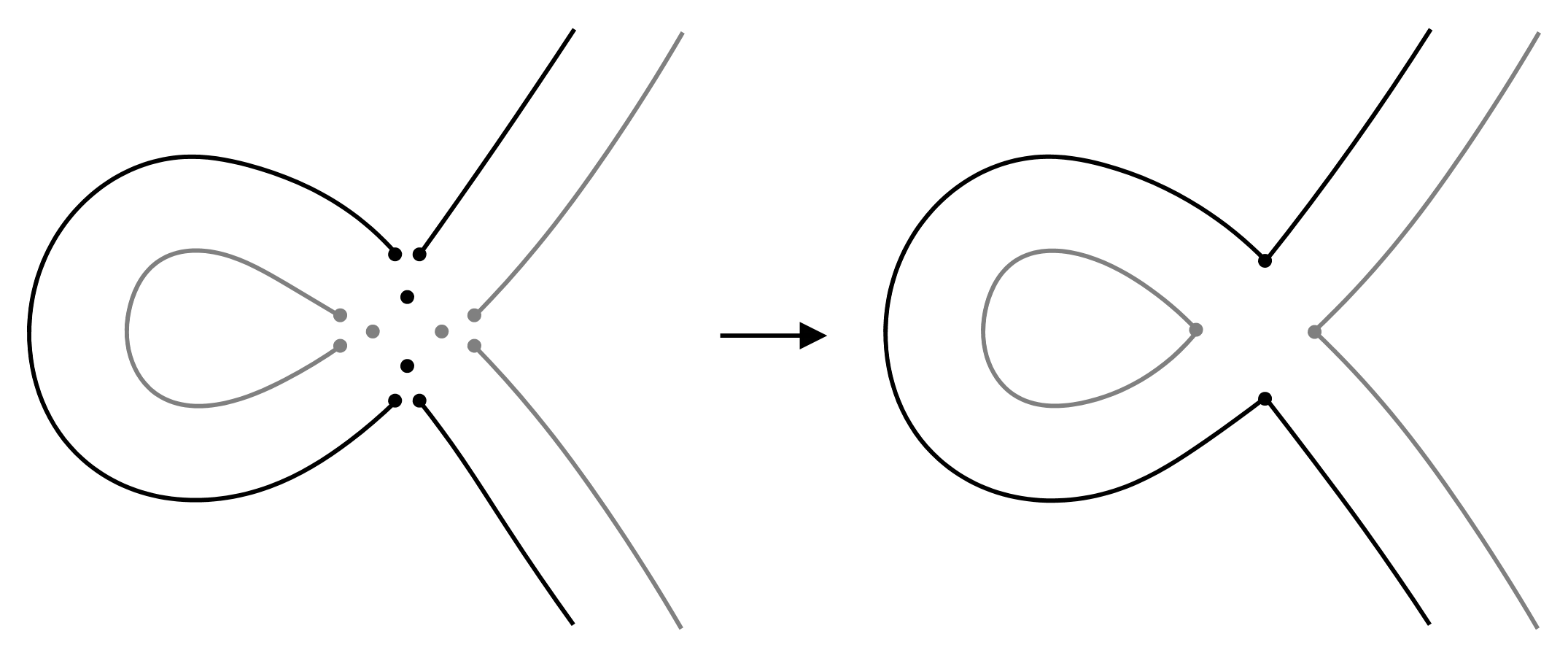}
	\end{minipage}\hfill
	\begin{minipage}[l]{0.37\textwidth}
		\caption{The gluing process for $\pi(x,y) = y^2 - x^3 - 3 x^2$. On the right hand side, 
	we see the initial spaces $\widehat{X}^+_0$ and $\widehat{X}^-_0$ for positive and negative strata, respectively. 
	The left hand side depicts $\BlX^+_0$ and $\BlX^-_0$.}
		\label{gluing}
	\end{minipage}		
  \end{figure}%
\end{example}

\begin{remark} \label{rem:gluing}
  Let us make a few remarks.
  \begin{enumerate}
		\item 
		  Of course, we could alternatively describe the gluing starting only with strata of maximal dimension. 
		  For example, $\BlX$ is the space obtained from
			\[
				\widehat{X}' = \bigsqcup_{\substack{S \in \SS \\ \codim(S) = 0}} \widehat{S} \times Q(S)
			\]
			via the following equivalence relation: Two points 
			\begin{gather*}
			  (p_1,e_1,q_1) \in T_1 \times E(T_1,S_1) \times Q(S_1) \subset \widehat{X}', \\
			  (p_2,e_2,q_2) \in T_2 \times E(T_2,S_2) \times Q(S_2) \subset \widehat{X}',
			\end{gather*}
			are declared equivalent if
			$T_1 = T_2=T$, $p_1 = p_2=p$ and $\mu_{T,S_1}(e_1, q_1) = \mu_{T,S_2}(e_2, q_2) \in Q(T)$. 
		\item 
		  The situation is even simpler if we additionally assume that $D'$ is \emph{strictly} normal crossing,
			that is, all its irreducible components are smooth submanifolds of $X$. Here, we use the adhoc definition
			of irreducible component as union of closures of maximal strata that are connected in codimension $1$. 
			In this case, for $T \subsetneq S$ we have $|E(T,S)| = 1$. Hence we have a canonical inclusion 
			$Q(S) \hookrightarrow Q(T)$ and $\widehat{S} = \overline{S} \subset X$.
			We then glue $\overline{S_1}$ and $\overline{S_2}$ along $T \subsetneq S_1, S_2$ for sheets 
			labelled by $q_1 \in Q(S_1)$ and $q_2 \in Q(S_2)$ if and only if $q_1 = q_2 \in Q(T)$. 
		\item 
			The orthant sets $Q(S)$ carry the structure of affine spaces with tangent space
		  $\Z_2^\HH$, where $\Z_2 = \Z/2\Z$ and $\HH$ is the set of hyperplanes in $T_p D$ for a points $p \in S$. 
			The maps $\mu(e,.)$ are affine maps. 
			The positive orthants $Q^+(S) \subset Q(S)$ form an affine hyperplane tangent to 
			the hyperplane of points in $\Z_2^\HH$ whose coordinates sum up to zero. 		
		\item 
		  Again, \autoref{thm:GluingSmooth} generalises to the case of weakly semi-stable degenerations
			from \autoref{rem:weaklySemiStable}. The only change is that now, according to the 
			three cases in \autoref{eq:weaklyTypes}, we have $|Q^+(S)| = 2^{k-1}$, $Q^+(S) = Q(S)$ or
			$Q^+(S) = \emptyset$, respectively. However, the orthant maps still respect positive orthants
			and the gluing can still be described as above. 
		  
	\end{enumerate}
\end{remark}

\begin{example} \label{ex:deformationconic}
  Using the notation from \autoref{ex:normalcone}, 
	let us consider the semi-stable degeneration $X$ given 
	by the deformation of $Y = \C\P^2$ to the normal cone 
	of a conic $Z \subset \C\P^2$.
	In this case the special fibre $X_0$ consists of the two components
	$B = \C\P^2$ and $C = \P(\OO_Z(4) \oplus \OO_Z) = \Sigma_4$ (the fourth Hirzebruch surface)
	glued	along $E = Z = \C\P^1$. 
	Note that $\R \Sigma_4$ is a topological $2$-torus. 
	Hence the stratification of $\R X_0$ consists of the strata
	$S_E = \R E = S^1$, the two connected components of $\R B \setminus \R E$ (a disc $S_B^1$ and a Möbius strip $S_B^2$) and $S_C = \R C \setminus \R E = 
	N_{\R Z} \R Y$
	(a cylinder).
	In contrast to \autoref{ex:doublecovering}, the normal bundles $N_{\R Z} \R Y \approx N_{\R E} \R B \approx N_{\R E} \R C$
	are (topologically) trivial and hence $\wt{S}^+_E \to S_E$ is the trivial (disconnected) covering map of degree $2$. In particular, 
	even though the stratification is not simply connected, there is no monodromy and we can directly apply \autoref{thm:GluingSmooth}.
	The stratified space $\BlX^+_0$ is depicted in \autoref{normalcone}. 
	\begin{figure}[b]
		\begin{minipage}[r]{0.45\textwidth}
			\input{pic/normalcone.TpX}
		\end{minipage}\hfill
		\begin{minipage}[l]{0.5\textwidth}
			\caption{The stratified positive special fibre $\BlX^+_0$ for the deformation of $\C\P^2$ to the normal cone of a conic $Z \subset \C\P^2$. 
			         We see three copies of the real conic $\R Z$. The inner and outer one together form $\wt{S}^+_E$. The middle one represents the zero section
							of $\R N_Z Y = N_{\R Z} \R Y$. In essence, the deformation replaces $\R Z$ by a closed tubular neighbourhood whose boundary is $\wt{S}^+_E$ and whose 
							interior is $N_{\R Z} \R Y$.}
			\label{normalcone}
		\end{minipage}		
	\end{figure}
	
	We see that this deformation can potentially be used to \enquote{embed}
	real algebraic curves in $\Sigma_4$ (in particular, those which do not intersect $E=Z$) inside $\C\P^2$. 
	For example, in \cite{Che-FourMCurves, Ore-NewMCurve} certain real algebraic curves in $\Sigma_4$ are used to construct new
	$M$-curves of degree $8$ in $\C\P^2$. 
	The transfer of these curves from $\Sigma_4$ to $\C\P^2$ is achieved via 
	the perturbation of the singularity of four maximally tangent conics
	which can also be interpreted/realized via 
	the above deformation to the normal cone, c.f.\ \cite[Section 7]{Shu-LowerDeformationsIsolated}
	and \cite[Chapter 3]{ST-PatchworkingSingularAlgebraic}. 
	In the latter papers, deformations to the normal cone are used systematically to transfer 
	a local deformation of singularities to a deformation of the whole curve. 
	In \cite[Proposition 4.7]{BDIM-RealAlgebraicCurves}, the deformation to the normal cone of a cubic curve 
	is used to construct a real algebraic curve $C \subset \C\P^2$
  of degree 12 such that $\R C$ consists of 45 isolated points. 
	(Of course, in all these examples additional deformation-theoretic arguments are needed 
	to prove that, in our notation, a given curve $C_0 \in N_Z Y \subset X_0$ can be extended to a transversal $C \subset X$
	in the sense of \autoref{thm:RelativeReal}). 
	
	Interestingly, 
	in contrast to \autoref{ex:doublecovering}, deformations to normal cones of curves of higher degree are not of Viro patchworking type, cf.\ \autoref{ex:linepatchworking}.
	In particular, the embedding steps mentioned above cannot be directly achieved via patchworking 
	(while patchworking is used for example in \cite{Che-FourMCurves, Ore-NewMCurve} to construct the curves in $\Sigma_4$ to start with). 
	My thanks go to Erwan Brugallé for sharing and explaining these ideas to me. 
\end{example}

\section{Toric degenerations} \label{toricdegenerations}

We now consider the particular case of toric degenerations (we use the term in the sense that $\pi$
is required to be a toric morphism). This example is not necessarily of interest 
by itself since it can be treated more directly, but it serves as an important ambient framework 
for applying our method in the context of tropicalisation as in 
\cite{Bru-EulerCharacteristicSignature} and our upcoming joint work \cite{RRS-RealPhaseStructuresb}.
We recall quickly the setup. More details can be found in \cite{Smi-TorusSchemesOver, Spe-TropicalLinearSpaces}. 

Let $\Sigma$ be a finite complete polyhedral subdivision of $\R^n$. We denote by $C(\Sigma)$ the fan 
in $\R^{n+1}$ obtained by taking the cone over $\Sigma \times \{1\}$ and completing it on level $0$
by the recession fan $\Omega = \Rec(\Sigma)$ of $\Sigma$. We call $\Sigma$ \emph{unimodular} if $C(\Sigma)$ is unimodular
(in particular, $\Sigma$ is $\Q$-rational in such case). 
We call $\Sigma$ \emph{strongly unimodular} if the last coordinate $x_{n+1}$ of the primitive generator
of any ray of $C(\Sigma)$ is either $0$ or $1$ (in particular, $\Sigma$ is $\Z$-rational in such case). 
This can be rephrased as follows: for every $\sigma \in \Sigma$ there exists $x_0 \in \Z^n$ and a part of a lattice
basis $e_1, \dots, e_k, f_1, \dots, f_l \in \Z^n$ such that 
\begin{align} 
  \sigma = S + C, & & S = \text{Conv}(x_0, x_0 + e_1, \dots, x_0 +e_k), & & C = \R_{\geq 0} \langle f_1, \dots, f_l\rangle. 
\end{align}
Note that for any unimodular $\Sigma$ there exists an integer $d > 0$ such that $d \Sigma$ is strongly unimodular. 

We denote by $X = \C\Sigma$ and $\R X = \R\Sigma$ the complex and real toric varieties associated to $C(\Sigma)$.
We denote by $\pi \colon \C\Sigma \to \C$ the canonical toric morphism corresponding to the
projection to the last coordinate $l \colon \R^{n+1} \to \R$.
Its generic fibre is equal to the $\C\Omega$, the toric variety associated to the 
recession fan. The special fibre $X_0$ is a union of torus orbits $\OO(\sigma)$ labelled
by the cells $\sigma \in \Sigma$ and naturally identified with the group of semigroup homomorphisms
\begin{align*} 
  \OO(\sigma) = \Hom(\sigma^\perp, \C^\times) && \text{with real locus} && \R\OO(\sigma) = \Hom(\sigma^\perp, \R^\times).
\end{align*}
The closure of such a complex or real torus orbit is denoted
by $X(\sigma)$ and $\R X(\sigma)$, respectively, and is equal to the complex or real toric variety 
associated to the star fan around $\sigma$ denoted by $\Star_\Sigma(\sigma)$.  
We denote by $T_\Z(\sigma)$ the tangent space of a rational polyhedron $\sigma \subset \R^n$ intersected with $\Z^n$ and set
for any field $F$ 
\[
  F(\sigma) := F^n / ( F \otimes T_\Z(\sigma)).
\]
Note that $\Star_\Sigma(\sigma) \subset \R(\sigma)$. 

If $\Sigma$ is unimodular, 
the variety	$X$ is smooth, the map $\pi$ is regular away from $0$ and $X_0 \subset X$ is a normal crossing divisor. 
In other words, in such case $\pi$ is weakly semi-stable. 

\begin{lemma} \label{lem:StronglyUnimodular}
  If $\Sigma$ is strongly unimodular, $\pi$ is a totally real semi-stable degeneration.
\end{lemma}

\begin{proof}
  This seems to be well-known, we sketch the argument showing that $\pi$ is normal crossing. 
	Given $\sigma \in \Sigma$, let
	\[
	  (x_0, 1), \dots, (x_k, 1), (f_1, 0), \dots, (f_l, 0)
	\]
	denote the primitive generators of the rays of $C(\sigma)$. Let $l_0, \dots, l_k$ linear forms on $\R^{n+1}$ defined over $\Z$ 
	such that $l_i((x_j, 1)) = \delta_{ij}$ and $l_i((f_j, 0))$ for all $i,j$. These linear forms correspond to monomial functions
	in a neighbourhood of $X(\sigma)$ generating the ideals of the $k$ irreducible components of $X_0$ containing $X(\sigma)$. In the same way, 
	the projection $l \colon \R^{n+1} \to \R$ corresponds to the function $\pi$. 
	But $l = l_1 + \dots + l_k$ on $C(\sigma)$, which proves that $\pi$ is normal crossing locally around $X(\sigma)$. 
\end{proof}

\begin{remark} 
  If $\Sigma$ is unimodular, we can write 
	\[
	  l = a_1 l_1 + \dots + a_k l_k 
	\]
	for positive integers $a_1, \dots, a_k$ which correspond to the integers mentioned in \autoref{rem:weaklySemiStable}.
	Note that since $l$ is primitive, at least one of the $a_i$ is odd and hence all local forms are of the first
	type as listed in \autoref{eq:weaklyTypes}. Hence the positive orthants behave as in the semi-stable case 
	and the following discussion also applies to the weakly semi-stable case. 
	In fact, the discussion in principle even applies to general subdivisions $\Sigma$, but we do not need this here. 
\end{remark}

We assume from now on that $\Sigma$ is strongly unimodular. 
The toric boundary divisors of $\C\Omega$ (that is, the rays of $\Omega$) give rise 
to a transversal collection of smooth divisors $X_1, \dots, X_m \subset X$.
We denote by $\SS = \SS(\R X, \R D')$ the stratification of $\R X$ 
associated to the normal crossing divisor
\[
  D' = X_0 \cup X_1 \cup \dots \cup X_m.
\]
By construction, $\SS$ is a refinement of $\SS(\R X, \R X_0)$. 
We denote by $\SS_0$ the subset of strata contained in $\R X_0$.  
These strata are exactly the connected components of the real torus orbits $\R \OO(\sigma) = \Hom(\sigma^\perp, \R^\times)$,
hence they are (open) orthants. 
Moreover, these orthants are naturally labelled by the elements of $\Z_2(\sigma)$.
In summary, the strata in $\SS_0$ are labelled by tuples $(\sigma, \epsilon)$
with $\sigma \in \Sigma$ and $\epsilon \in \Z_2(\sigma)$. 

The semigroup homomorphisms of absolute value 
$| \cdot | \colon \R^\times \to \R_{> 0}$ and logarithm 
$\log \colon \R_{> 0} \to \R$ induce maps 
$\Hom(\sigma^\perp, \R^\times) \to \Hom(\sigma^\perp, \R_{> 0})$ and
$\Hom(\sigma^\perp, \R_{> 0}) \to \Hom(\sigma^\perp, \R)$ which we denote by
the same symbols by abuse of notation. 
Note that here $\R$ is the additive (semi-)group and hence $\Hom(\sigma^\perp, \R)$ 
is the group of ordinary linear homomorphisms
For each stratum $S$ labelled by $\sigma$, this yields canonical homeomorphisms
\[
  S \stackrel{| \cdot |}{\approx} \Hom(\sigma^\perp, \R_{> 0}) \stackrel{\log}{\approx} \Hom(\sigma^\perp, \R) = \R(\sigma).
\]

Note that the normal crossing divisor $D'$ is strictly normal crossing. 
In particular, the endpoint modification $\widehat{S}$
of a stratum $S \in \SS_0$ is equal to the closure $\overline{S}$ in $\R X_0$, see \autoref{rem:gluing}. 
In fact, if $S$ is labelled by $(\sigma, \epsilon)$, it is sufficient to take the closure in
$\R X(\sigma)$, hence $\widehat{S}$ is a closed orthant of $\R X(\sigma)$
homeomorphic to the positive closed orthant. The latter can be described
via the usual semigroup formalism for toric varieties, but using 
the (multiplicative) semigroup $\R_{\geq 0}$.
This leads as naturally to tropical toric varieties which are defined
using again the same formalism, but this time  using the additive semigroup
$\T = \R \cup \{-\infty\}$, see
\cite[Section 3]{Pay-AnalytificationIsLimit} and \cite[Chapter 5]{MR-TropicalGeometry}.
Note that the extended logarithm map $\log \colon \R_{\geq 0} \to \T = \R \cup \{-\infty\}$
(setting $\log(0) = -\infty$) is a semigroup isomorphism and a homeomorphism. 
It follows that the tropical variety associated to a fan is canonically homeomorphic
to any closed orthant of the corresponding real toric variety. 
Applied to our situation, the statement is that $\widehat{S} = \overline{S}$ 
is canonically homeomorphic to $\T X(\sigma)$, the 
tropical toric variety associated to $\Star_\Sigma(\sigma)$ with open dense 
(tropical) torus $\R(\sigma)$. 

\begin{corollary} \label{cor:toricgluing}
  The space $\BlX^+_0$ is obtained from 
	\begin{equation} \label{eq:startingspaceToric} 
	  \bigsqcup_{\sigma \in \Sigma} \T X(\sigma) \times \Z_2(\Rec(\sigma))
	\end{equation}
	via gluing generated by the following:
	Given pairs $(\sigma, \epsilon)$ and $(\sigma', \epsilon)$, we identify
	$\T X(\sigma) \times \{\epsilon\}$ with the copy contained in 
	$\T X(\sigma') \times \{\epsilon'\}$ if 
	\begin{align*} 
		\sigma' \subset \sigma && \text{and} && \epsilon' \mapsto \epsilon && \text{under} && \Z_2(\Rec(\sigma')) \to \Z_2(\Rec(\sigma)). 
	\end{align*}
\end{corollary}

\begin{proof}
  The statement is the result of applying \autoref{thm:GluingSmooth} to our particular case. 
	In order to see this, it is enough to understand how the initial space from 
	\autoref{eq:startingspaceToric} is identified with the initial space from
	\autoref{eq:startingspaceSmooth}. To this end, let 
	us give a more intrinsic description of the orthants
	$Q(S)$ of a stratum $S$ labelled by $(\sigma, \epsilon)$. 
	It is useful to first study the set of orthants $\wtQ(S)$ with respect to the divisor $D'$ (instead
	of $X_0$ for $Q(S)$). 
	In fact, note that the maximal strata of $\SS$ 
	(equivalently, the orthants of $\R X = \R\Sigma$)
	are globally labelled by $\Z_2^{n+1}$. 
	Given $p \in S$, the tangent cone $T_p D' \subset T_p X$ consists of $k'$ hyperplanes
	which are naturally labelled by the rays of $C(\sigma) \subset \R^{n+1}$ (the cone over
	$\sigma \times \{1\}$). The action on $\wtQ(S)$ given by crossing a hyperplane (denoted by $\Z_2^\HH$ in \autoref{rem:gluing})
	under this identifications corresponds to the additive action of the the subspace $\Z_2 \otimes T_\Z(C(\sigma)) \subset \Z_2^{n+1}$. 
	Hence $\wtQ(S)$ is an affine subspace in $\Z_2^{n+1}$ with tangent space $\Z_2 \otimes T_\Z(C(\sigma))$.
	In fact, we can check easily that 
	\[
		\wtQ(S) = q_1^{-1}(\epsilon)  \; \subset \Z_2^{n+1}
	\]
	where $q$ denotes the quotient map
	\[
		q_1 \colon \Z_2^{n+1} \to \Z_2(C(\sigma)) = \Z_2(\sigma).
	\]
	Now, in order to pass from $\wtQ(S)$ to $Q(S)$, we just have to quotient by
	the subaction corresponding to the divisors among $X_1, \dots, X_m$ that contain $S$. 
	These correspond to the rays in $C(\Sigma) \cap (\R^n \times \{0\}) = \Rec(\sigma) \times \{0\}$, 
	where $\Rec(\sigma)$ denotes the recession cone of $\sigma$.
	In summary, 
	\begin{align*} 
		Q(S) = q_2^{-1}(\epsilon) && \text{with} && q_2 \colon \Z_2(\Rec(\sigma) \times\{0\}) \to \Z_2(C(\sigma)) = \Z_2(\sigma).
	\end{align*}
	Finally, let us restrict to positive orthants. Clearly, the set of positive orthants in $\Z_2^{n+1}$ 
	is the kernel of the map $\Z_2^{n+1} \to \Z_2$ given by summing all coordinates.
	By projecting to the first $n$ coordinates, we can hence identify positive orthants with $\Z_2^n$. 
	With this last identification, we obtain
	\begin{align*} 
		Q^+(S) = q^{-1}(\epsilon) && \text{with} && q \colon \Z_2(\Rec(\sigma)) \to \Z_2(\sigma).
	\end{align*}
	Summing over all strata labelled by a fixed $\sigma$ of sedentarity $s$, we obtain
	\[
		\bigsqcup_{\substack{S \in \SS \\ S \subset \OO(\sigma)}} \overline{S} \times Q^+(S) = \T X(\sigma) \times \Z_2(\Rec(\sigma)).
	\]
	This provides the identification between the two initial spaces in question.
	It is straightforward to check that the gluing described in \autoref{thm:GluingSmooth}
	coincides with the gluing described here. 
\end{proof}

\begin{example} \label{ex:linepatchworking}
  The degeneration from \autoref{ex:doublecovering} is in fact a toric degeneration induced 
	by the subdivision $\Sigma$ depicted in \autoref{toricP2}. Moreover, this subdivision
	is dual to a convex subdivision of a lattice polytope as used in Viro's patchworking, 
	cf.\ \autoref{ex:patchworking}.
	Note that in contrast to the stratification $\SS(\R X, \R X_0)$ discussed in \autoref{ex:doublecovering}, 
	the stratification $\SS$ that takes into account the toric divisors is simply connected and the gluing 
	procedure is described by \autoref{cor:toricgluing}. %
	\begin{figure}[ht]
		\begin{minipage}[r]{0.47\textwidth}
			\includegraphics[width=\textwidth]{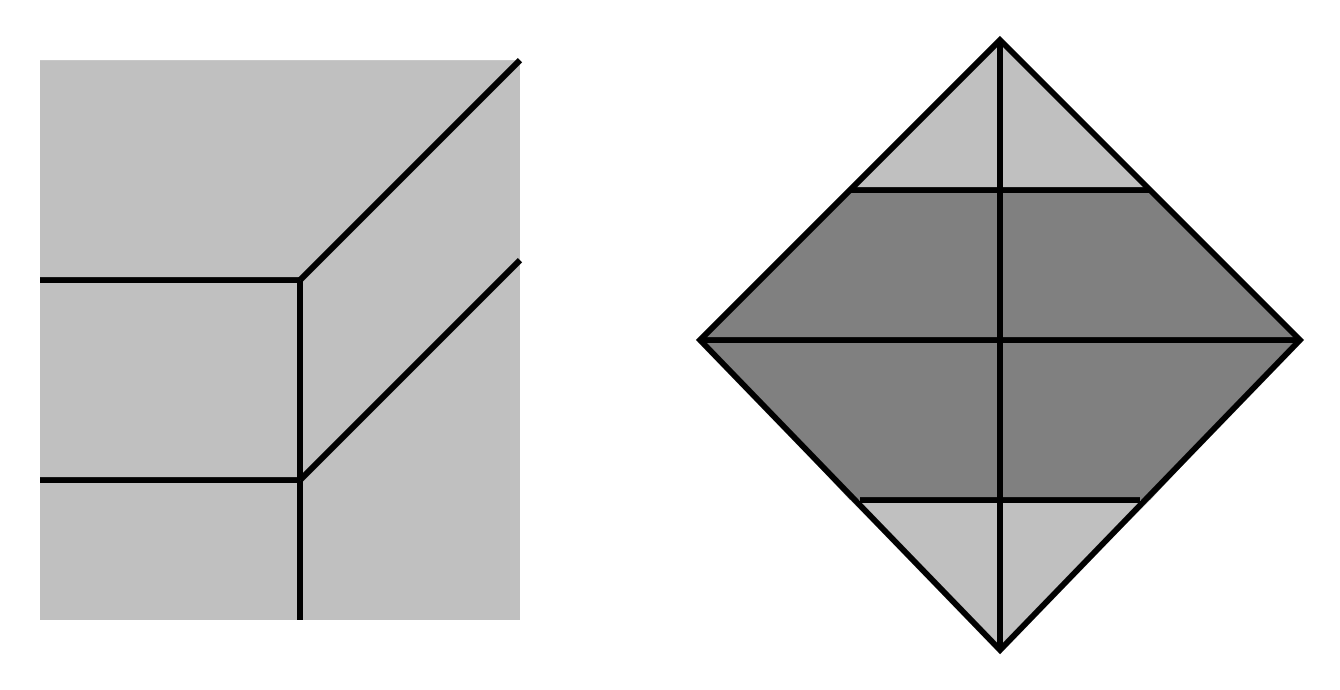}
		\end{minipage}\hfill
		\begin{minipage}[l]{0.5\textwidth}
			\caption{The left hand side depicts the subdivision $\Sigma$ that produces the degeneration from 
							 \autoref{ex:doublecovering}. The right hand side represents the gluing of $\BlX^+_0 \approx \R\P^2$ via 
							 four copies of $\TP^2$ and four copies of $\Bl_p \TP^2$ 
							 (one should additionally glue opposite boundary points).}
			\label{toricP2}
		\end{minipage}		
	\end{figure}
\end{example}

As final ingredient of our discussion, we now add a real submanifold $Y \subset X$ 
such that $Y, X_1, \dots, X_m$ is a collection
of transversal submanifolds. We call such a $Y$ \emph{torically transversal}. 
For any stratum $S \in \SS_0$ labelled by $(\sigma, \epsilon)$ (we do \emph{not} change the stratification), we denote 
\begin{align*} 
  Y(\sigma, \epsilon) = Y \cap S, \\
	\overline{Y}(\sigma, \epsilon) = Y \cap \overline{S}. 
\end{align*}
Using the identifications from above, we consider these sets as subsets of $\R(\sigma)$ and $\T X(\sigma)$, 
respectively. 
We state for convenience the following summary which is just the combination of \autoref{thm:RelativeReal} and \autoref{cor:toricgluing}
applied to our case. 

\begin{corollary} \label{cor:tropicalisation}
  The topological pair $\R\wt{Y}^+_0 \subset \R\BlX^+_0$ is obtained via gluing as described in \autoref{cor:toricgluing} from the topological pair
	\[
		\bigsqcup_{\substack{\sigma \in \Sigma \\ \epsilon \in \Z_2(\Rec(\sigma))}} \overline{Y}(\sigma, [\epsilon]) \times \{\epsilon\} \;\;\; \subset \;\;\; \bigsqcup_{\sigma \in \Sigma} \T X(\sigma) \times \Z_2(\Rec(\sigma)).
	\]
	Here, $[\epsilon]$ denotes the projection to $\Z_2(\sigma)$.
	Moreover, the topological pair $\R\wt{Y}^+_0 \subset \R\BlX^+_0$ 
	is homeomorphic to $\R Y_t \subset \R X_t$ for all $t \in (0,1)$
	and these homeomorphisms can be chosen such that they
	respect the intersections with the divisors $X_1, \dots, X_m$. 
\end{corollary}

\begin{example} \label{ex:patchworking}
  Viro's patchworking method \cite{Vir-CurvesDegree7, Vir-PatchworkingRealAlgebraic} 
	takes as input a Laurent polynomial in $n$ variables $F_0(x) = \sum_{I \in \Z^n} a_I x^I$
	and a convex subdivision $\Gamma$ of the Newton polytope $\Delta$ of $F_0$ 
	(induced by a piecewise linear convex function $\nu \colon \Delta \to \R$; we may 
	assume that $\nu$ takes integer values on integer points). 
	The associated \emph{Viro polynomial} is 
	\[
	  F(x,t) = \sum_{I \in \Delta \cap \Z^n} t^{\nu(I)} a_I x^I \quad \quad \in \R[x_1^\pm, \dots, x_n^\pm, t^\pm]. 
	\]
	For $t \neq 0$, we denote by $V_t$ the closure of $\{F(x,t) = 0\}$ in $\C\Delta$, the toric variety associated to $\Delta$. %
  \begin{figure}[t]
	\begin{minipage}[r]{0.6\textwidth}
		\includegraphics[width=\textwidth]{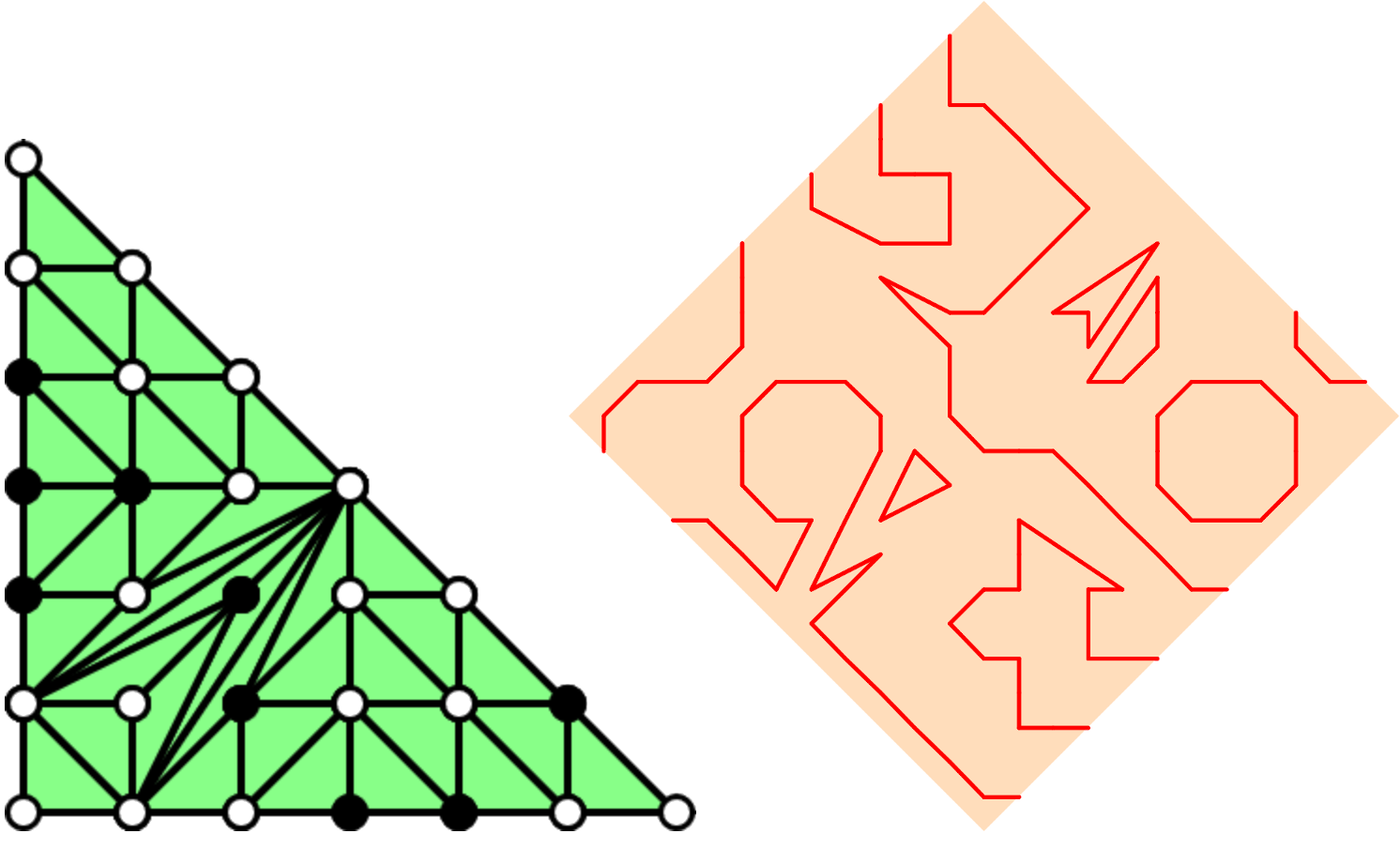}
	\end{minipage}\hfill
	\begin{minipage}[l]{0.37\textwidth}
		\caption{A unimodular triangulation of the the simplex of size $6$ and the associated patchwork. According to Viro's  patchwork theorem, it represents 
		the topological type of a real algebraic curve of degree $6$. The pictures were created using 
		\cite{HRR-CombinatorialPatchworkingTool}.}
		\label{patchworking}
	\end{minipage}		
  \end{figure}%
	Let us furthermore assume that $(F_0, \Gamma)$ is \emph{non-degenerate}, that is, for any $\sigma \in \Gamma$ 
	the truncation $F_0|_\sigma$ defines a non-singular hypersurface in $(\C^*)^n$.
	Under this assumption, Viro's patchworking theorem \cite[Theorem 1.7.A]{Vir-PatchworkingRealAlgebraic} provides a description of
	$\R V_t \subset \R\Delta$ for small values of $t$ in terms of gluing \emph{patches} given by
	the truncations $F_0|_\sigma$, $\sigma \in \Gamma$. 
	An example of the particular case of combinatorial unimodular patchworking (that is, $\Gamma$ is a unimodular triangulation of $\Delta$)
	is depicted in \autoref{patchworking}. 
	The connection to our setup is as follows:
	Let $\overline{\Delta}$ the upper graph of the convex function associated to $\nu$. 
	Let $\widetilde{\Delta}$ be a \emph{desingularisation} $\overline{\Delta}$ of $F$, that is, the normal fan of 
	$\widetilde{\Delta}$ is unimodular and refines the normal fan 
	of $\overline{\Delta}$.
	Then the projection to $t$ induces a map $\pi \colon \C \widetilde{\Delta} \to \C$
	which is a toric totally real weakly semi-stable degeneration 
	with generic fibre $\C\Delta$ as explained above. 
	Moreover the closure $V$ of $\{F(x,t) = 0\}$ in $\C \widetilde{\Delta}$ is 
	a real submanifold which is torically transversal.
	This is essentially the statement of \cite[Theorem 2.5.C]{Vir-PatchworkingRealAlgebraic},
	even though the terminology of semi-stable (and toric) degenerations is not used explicitly there. 
	Hence our methods recover\fshyp{}reformulate parts of 
	Viro's patchworking method
	(see \cite[Theorem 3.3.A]{Vir-PatchworkingRealAlgebraic} for connecting 
	the gluing from \autoref{cor:tropicalisation} with respect to the desingularisation $\Delta$ 
	to the patchwork described in \cite[Section 1.5]{Vir-PatchworkingRealAlgebraic}).
	We hope to generalise this to tropicalisations of arbitrary codimension in
	our upcoming joint work \cite{RRS-RealPhaseStructuresb} with Arthur Renaudineau and Kris Shaw.
\end{example}

\printbibliography

\vfill
\section*{Contact}

    Johannes Rau \\ Departamento de Matemáticas \\ Universidad de los Andes \\
		KR 1 No 18 A-10, BL H \\ Bogotá, Colombia \\ \href{mailto:j.rau@uniandes.edu.co}{j.rau AT uniandes.edu.co}

\end {document}